\documentclass[11pt]{amsart}

\usepackage[T1]{fontenc}
\usepackage{amssymb,amsmath,amsfonts,amsthm,mathtools}
\usepackage{microtype,booktabs,array,tabularx,longtable,enumitem}
\usepackage{placeins}
\usepackage[font=small]{caption}
\usepackage{xcolor}
\usepackage{hyperref}
\usepackage{cleveref}

\definecolor{darkblue}{rgb}{0.0,0,0.7}
\hypersetup{
 colorlinks=true,
 citecolor=darkblue,
 linkcolor=darkblue,
 urlcolor=darkblue,
 pdftitle={Strongly explicit Ramanujan 2-lift towers via Hilbert modular forms},
 pdfsubject={Binary Ramanujan towers for primes congruent to 1 modulo 17},
 pdfauthor={Olga Kharlampovich; Igor Pak; Alina Vdovina},
 pdfkeywords={Ramanujan graphs, 2-lifts, strong explicitness, Hilbert modular forms, Jacquet-Langlands correspondence, Hecke operators}}

\makeatletter
\def\th@plain{%
  \thm@notefont{}%
  \itshape
}
\def\th@definition{%
  \thm@notefont{}%
  \normalfont
}
\makeatother

\setlist[enumerate]{label=\textnormal{(\roman*)},ref=\roman*}

\newtheorem{theorem}{Theorem}[section]
\newtheorem{proposition}[theorem]{Proposition}
\newtheorem{lemma}[theorem]{Lemma}
\newtheorem{corollary}[theorem]{Corollary}
\theoremstyle{definition}

\newtheorem{remark}[theorem]{Remark}
\newcommand{\Q}{\mathbb Q}
\newcommand{\Z}{\mathbb Z}
\newcommand{\C}{\mathbb C}
\newcommand{\F}{\mathbb F}
\newcommand{\OF}{\mathcal O_F}
\newcommand{\OO}{\mathcal O}
\newcommand{\EE}{\mathcal E}
\newcommand{\levelprime}{Q}
\newcommand{\otherlevelprime}{Q'}
\newcommand{\adjprime}{P}
\newcommand{\localprime}{R}
\newcommand{\qi}{\mathbf{i}}
\newcommand{\qj}{\mathbf{j}}
\newcommand{\qk}{\mathbf{k}}
\newcommand{\splitmap}{\psi}
\newcommand{\mobius}{\varphi}
\newcommand{\fundunit}{\varepsilon_0}
\DeclareMathOperator{\red}{red}
\newcommand{\one}{\mathbf 1}
\DeclareMathOperator{\nrd}{nrd}
\DeclareMathOperator{\trd}{trd}
\DeclareMathOperator{\Cl}{Cl}
\DeclareMathOperator{\GL}{GL}
\DeclareMathOperator{\PGL}{PGL}
\DeclareMathOperator{\Spec}{Spec}
\DeclareMathOperator{\Tr}{Tr}
\newcommand{\IdealNorm}[1]{\mathrm{N}\!\left(#1\right)}
\newcommand{\FieldNorm}[1]{\mathrm{N}_{F/\Q}\!\left(#1\right)}
\DeclareMathOperator{\diag}{diag}
\newcolumntype{Y}{>{\raggedright\arraybackslash}X}
\newcolumntype{P}[1]{>{\raggedright\arraybackslash}p{#1}}
\numberwithin{equation}{section}

\crefname{theorem}{Theorem}{Theorems}
\crefname{lemma}{Lemma}{Lemmas}
\crefname{proposition}{Proposition}{Propositions}
\crefname{definition}{Definition}{Definitions}
\crefname{corollary}{Corollary}{Corollaries}
\crefname{section}{Section}{Sections}
\crefname{subsection}{Section}{Sections}
\crefname{table}{Table}{Tables}
\crefname{equation}{equation}{equations}
\newcommand{\arxivref}[1]{\href{https://arxiv.org/abs/#1}{\texttt{arXiv:#1}}}
\makeatletter
\renewcommand{\@biblabel}[1]{[#1]\hfill}
\makeatother

\title[Strongly explicit Ramanujan $2$-lift towers]{Strongly explicit Ramanujan $2$-lift towers \\ via Hilbert modular forms}
\date{\today}

\author{Olga Kharlampovich}
\address[Olga Kharlampovich]{CUNY Hunter College and Graduate Center}
\email{\texttt{okharlam@hunter.cuny.edu}}

\author{Igor Pak}
\address[Igor Pak]{University of California, Los Angeles}
\email{\texttt{pak@math.ucla.edu}}

\author{Alina Vdovina}
\address[Alina Vdovina]{CUNY City College and Graduate Center}
\email{\texttt{avdovina@ccny.cuny.edu}}

\begin{document}

\begin{abstract}
Expander graphs have been central objects of study for over fifty years.
\emph{Ramanujan graphs} attain the asymptotically optimal nontrivial adjacency-eigenvalue
threshold for infinite fixed-degree families.  The original arithmetic constructions by
Lubotzky--Phillips--Sarnak (1988) and Margulis (1988), and the iterative $2$-lift
method of Bilu--Linial (2006) provide different approaches to Ramanujan graph
constructions.  Motivated by work on Mealy automata, we combine both approaches to
give the following new arithmetic construction.

For every prime $p\equiv 1\pmod {17}$, we construct a infinite tower of $(p+1)$-regular non-bipartite
Ramanujan graphs on $2^n$ vertices, with label-preserving $2$-lifts between consecutive levels.  These construction is \emph{strongly explicit}: for fixed $p$, a neighbor or lift sign is computable in $O_p(n^3)$ bit operations. The construction uses fixed algebraic $2$-adic transformations, specified by finite data computable in preprocessing depending only on $p$.

The main number-theoretic tool is Blasius's Ramanujan--Petersson theorem for Hilbert modular forms, applied through the Jacquet--Langlands correspondence to obtain the sharp bound $|\lambda|\leq2\sqrt p$ for every nontrivial adjacency eigenvalue of our $(p+1)$-regular graphs. The transformations preserve $2$-adic distance, which gives the binary covers and an explicit formula for their signings.
\end{abstract}

\maketitle

\section{Introduction}\label{sec:intro}


\subsection{Foreword}
The story of expanders is both extensive and extremely surprising, as their applications
are crucial in both algorithms and computational complexity, while the tools involved
are quite technical and substantially utilize many areas of mathematics \cite{HLW06,Lub12,Vad12}.

Three developments are particularly relevant here: the explicit expanders of Margulis \cite{Mar73}; the arithmetic Ramanujan graphs of Lubotzky--Phillips--Sarnak \cite{LPS88} and Margulis \cite{Mar88}; and the iterative $2$-lift approach of Bilu--Linial \cite{BL06}, further developed by Marcus--Spielman--Srivastava \cite{MSS15}.
These constructions have prompted many generalizations, often motivated by applications that require different combinations of spectral quality, graph size, and computational access. The literature overview below separates these requirements.



The goal of this paper is to obtain a strongly explicit construction of $2$-lifts
which have Ramanujan property, thus combining all three properties above.  This is
especially difficult since there is no ``room to recover'' in this situation --- once
an eigenvalue is outside of the Ramanujan interval, it remains so under $2$-lifts.
Furthermore, being strongly explicit and Ramanujan, we are essentially forced to
use number-theoretic tools, as both combinatorial and probabilistic approaches are
less powerful and typically give either $\varepsilon$-Ramanujan or probabilistic constructions
(see below).

The story of our result is also quite interesting.  Originally, the
authors made an effort to understand expansion properties of Aleshin and Bellaterra
Schreier graphs \cite{KV26a,KV26b,MP17}.  The construction in this paper grew
out of this work.  Although we retain the arithmetic model underlying the
cubic Bellaterra Ramanujan tower and extend its analysis, the resulting
graphs are no longer obtained as Schreier of automatic groups.




\subsection{Main result}
For introductions to Ramanujan graphs and the lifting method,
see Kowalski \cite[\S4.2]{Kow19} and Wigderson \cite[\S8.7]{Wig19}, respectively.
For the Marcus--Spielman--Srivastava method, see an accessible exposition
in \cite[\S13.3]{Wig19}.

When a $d$-regular expander $G=(V,E)$ has $|V|=2^n$ vertices,
a \emph{strongly explicit} construction computes $i$-th neighbor of a vertex $v\in V$,
where $1\le i \le d$, in time polynomial in~$n$, without constructing the entire graph.
We refer to Section~\ref{s:notation} for a more precise definition.


We give a strongly explicit construction of a \emph{compatible binary tower} $\{G_n\}$, in which $G_{n+1}$ on $(n+1)$-bit strings is a two-sheeted cover of $G_n$ on $n$-bit strings. Each base edge has a \emph{lift bit} $\beta\in\{0,1\}$, specifying the straight or crossed matching between its two fibers. The corresponding \emph{sign} is defined as $(-1)^\beta\in\{\pm 1\}$.

A connected $d$-regular non-bipartite graph with adjacency matrix $A$ is called \emph{Ramanujan}, if every eigenvalue $\lambda<d$ has absolute value $|\lambda|\le 2\sqrt{d-1}$.  By the Alon--Boppana theorem, this is the asymptotically optimal bound for infinite fixed-degree families, see e.g.\ \cite{HLW06}. 

\begin{theorem}[{\rm Strongly explicit Ramanujan towers}{}]\label[theorem]{thm:main}
For every prime $p\equiv1\pmod {17}$, there is a deterministically specified family
$\{\Gamma_{n,p}\}_{n\geq 1}$ of connected graphs with the following properties:
\begin{enumerate}[label=\textup{(\roman*)}]
\item $V(\Gamma_{n,p})=\Z/2^n\Z$, and the degree of every vertex is equal to \thinspace $d:=p+1$.
\item For every $n\geq1$, every adjacency eigenvalue $\lambda<d$ satisfies
\begin{equation*}
 |\lambda|\, \leq \, 2\sqrt p\,.
\end{equation*}
Thus every $\Gamma_{n,p}$ is a connected, non-bipartite $d$-regular Ramanujan graph.
\item Reduction modulo $2^n$ is a label-preserving two-sheeted cover
\[
 \Gamma_{n+1,p}\longrightarrow\Gamma_{n,p}\,.
\]
The associated straight/crossed signing is deterministically and locally computable.
\item For fixed $p$, a neighbor or a lift bit can be computed in $O_p(n^3)$ bit operations.
\end{enumerate}
\end{theorem}

In other words, we have a simultaneous, effective description of an \emph{entire binary tower}:
exact order $2^n$ at every level, exact two-sided spectral bounds, and locally computable neighbors
and lift signs.
Dirichlet's theorem gives infinitely many primes $p\equiv1\pmod {17}$, see e.g.\ \cite{Dav00}.
The first three are $103,137,239$; see e.g.\ \cite[\href{http://oeis.org/A129484}{A129484}]{OEIS}. We conclude:

\begin{corollary}[{\rm Infinitely many degrees}{}]\label[corollary]{cor:degrees}
There are infinitely many degrees $d$  for which every binary level has a strongly explicit Ramanujan graph and every successive map is a $2$-lift.
\end{corollary}

Finally, we note that graphs in the theorem are labeled multigraphs: parallel edges, loops,
and half-loops are permitted, with the conventions specified in Section~\ref{s:notation}.

\subsection{Prior work}\label{sec:literature}

We recall several strands of work that motivate the construction: arithmetic Ramanujan graphs, graph coverings and lifts, derandomization, and strongly explicit expanders. The notions and parameter regimes differ across these literatures, so the summary below is intended mainly to locate the ingredients used later in the paper.

\smallskip\noindent\emph{Coverings and the tower viewpoint.}
Graph coverings have a classical permutation-voltage description, in which each directed edge carries a permutation and the reversed edge carries its inverse; see e.g.\ Gross--Tucker \cite{GT87}. Characteristic polynomials of graph coverings were studied by Mizuno--Sato \cite{MS95}. Nested finite quotients of groups give another natural source of compatible graph families, closely related to property~$(\tau)$, see   \cite{Lub94,Lub12}.

\smallskip\noindent\emph{Arithmetic Ramanujan graphs.}
Lubotzky--Phillips--Sarnak \cite{LPS88} and Margulis \cite{Mar88} gave the first arithmetic Ramanujan families. Pizer \cite{Piz90} subsequently constructed Ramanujan multigraphs directly from Hecke operators on spaces of theta series, a particularly close antecedent to the Hecke-operator viewpoint used below. Chiu \cite{Chi92} obtained a cubic construction, and Morgenstern \cite{Mor94} constructed strongly explicit Ramanujan graphs of degree $q+1$ for every prime power~$q$. Alon \cite{Al21} discusses strong explicitness and the range of attainable vertex sizes in these and related constructions.

\smallskip\noindent\emph{Random lifts and the universal-cover benchmark.}
Amit--Linial \cite{AL06} analyzed expansion in random lifts (the conference version appeared in 2001). Soon after, Friedman \cite{Fri03} studied the new spectrum of random lifts relative to their universal cover. His later random-regular-graph theorem \cite{Fri08} is the starting point for Bordenave's new proof and extension to random lifts of a fixed base graph \cite{Bor20}. Linial--Puder \cite{LP10} and Puder \cite{Pud15} developed further spectral estimates using word maps and counting methods. Agarwal--Kolla--Madan \cite{AKM13} study small and shift lifts, and Agarwal--Chandrasekaran--Kolla--Madan \cite{ACKM19} develop a group-based lifting framework. Generalized Alon--Boppana bounds for finite quotients of a fixed universal cover are given by Friedman--Tillich \cite{FT05}.

\smallskip\noindent\emph{Derandomizing lifts and Ramanujan coverings.}
Bilu--Linial \cite{BL06} introduced the iterative $2$-lift approach and gave deterministic polynomial-time constructions with nontrivial eigenvalues $O\big(\sqrt d\,(\log d)^{3/2}\big)$. In \cite{MSS15}, Marcus--Spielman--Srivastava proved that every bipartite $d$-regular Ramanujan graph has a Ramanujan $2$-lift.  In \cite{LPV14}, Liu--Peyerimhoff--Vdovina proved analogous existence result for $3$-lifts by using Hermitian metrices, and Chandrasekaran--Velingker \cite{CV17} proved analogous existence results for $4$-lifts. Soon after, Marcus--Spielman--Srivastava all-sizes theorem in \cite{MSS18} used a different interlacing family.  In \cite{Coh16}, Cohen gave a deterministic polynomial-time implementation producing bipartite Ramanujan multigraphs.

In \cite{HPS18}, Hall--Puder--Sawin showed that every finite connected graph admits a one-sided Ramanujan $r$-cover for every $r\geq1$, with a two-sided bound when the base graph is bipartite. Iterating Ramanujan $2$-lifts already gives existential binary towers. Their discussion of composition and wreath products \cite[\S5.4]{HPS18} also connects to Makelov's tower constructions \cite{Mak15}.

\smallskip\noindent\emph{Recent prescribed-graph signing algorithms.}
Most recently, Jadbabaie--Saberi--Sra \cite[Thm~6.1 and \S6.6]{JSS26} give a Las Vegas recursive-repair algorithm for sharp signings of finite simple bipartite graphs, with expected polynomial bit complexity. Xu \cite{Xu26} builds on this framework and gives a deterministic worst-case polynomial-time algorithm for the same sharp bipartite bound. These are whole-graph algorithms for a prescribed finite input graph.
{\small
\begin{table}[!htbp]
\centering
\renewcommand{\arraystretch}{1.12}
\begin{tabularx}{\textwidth}{@{}P{1.15in}P{1.17in}P{1.35in}Y@{}}
\toprule
Construction & Spectral target & Algorithmic guarantee & Covering structure \\
\midrule
LPS; Morgenstern & Exact & Local query: $\operatorname{polylog}(N)$ & Standard formulations do not give iterated $2$-lifts. \\
\addlinespace
Bilu--Linial & $O\bigl(\sqrt d\,(\log d)^{3/2}\bigr)$ & Whole graph: $\operatorname{poly}(N)$ & Deterministic iterated $2$-lifts. \\
\addlinespace
MSS I; HPS & Exact if bipartite & Existence & Binary lifts; HPS allows every covering degree, with a general one-sided bound. \\
\addlinespace
MSS IV; Cohen & Exact & Whole graph: $\operatorname{poly}(N)$ & All admissible sizes; not a prescribed-base local signing algorithm. \\
\addlinespace
JSS; Xu & Exact, bipartite & Whole graph: $\operatorname{poly}(N)$ & Prescribed-graph signings; randomized repair, then deterministic repair. \\
\addlinespace
MOP & $2\sqrt{d-1}+\varepsilon$ & Whole graph: $\operatorname{poly}(N)$ & Derandomized lifts for every fixed degree. \\
\addlinespace
O'Donnell--Wu & $\varepsilon$-close to $\Spec(X)$ & Whole graph: $\operatorname{poly}(N)$ & Repeated $2$-lifts; a separate probabilistically strongly explicit variant. \\
\addlinespace
RSW; Makelov & Expanding; almost-Ramanujan & Local query: $\operatorname{polylog}(N)$ & Compatible towers; exact binary towers also treated existentially by Makelov. \\
\addlinespace
BGV; KV & Exact & Local query: $O(\log N)$ via fixed automata & Respectively $q$-fold towers for odd prime powers $q$, and a cubic binary tower. \\
\addlinespace
This paper & Exact & Local query: $O_p((\log N)^3)$ & Every size $N=2^n$; locally computed binary lift signs; $d=p+1$, $p$ prime. \\
\addlinespace
\bottomrule
\end{tabularx}
\caption{Selected guarantees for graphs on $N$ vertices. Degree and approximation parameters are fixed in the runtime bounds. ``Whole graph'' means a construction-time bound, not a local-query bound.  JSS gives an expected bound and Xu a worst-case bound \cite{JSS26,Xu26}. Graph categories and degree ranges differ. The almost-Ramanujan variants have degree-asymptotic guarantees, not an arbitrary additive error at each fixed degree.}
\label{tab:comparison}
\end{table}
}

\smallskip\noindent\emph{Near-Ramanujan graphs and combinatorial explicitness.}
In \cite{MOP20}, Mohanty--O'Donnell--Paredes gave deterministic polynomial-time constructions with bound $2\sqrt{d-1}+\varepsilon$ for every fixed $d\geq3$ and $\varepsilon>0$, using derandomized lifts. Jeronimo--Mittal--O'Donnell--Paredes--Tulsiani \cite{JMOPT22} extended closely related methods to explicit abelian lifts, with applications to classical and quantum LDPC codes. Alon \cite{Al21} gave strongly explicit constructions with bound $(2+o_d(1))\sqrt d$ for every degree. In \cite[\S1.3,~$\S$10.1]{OW20}, O'Donnell--Wu extended the lift approach to matrix-polynomial lifts, obtaining finite graphs whose nontrivial spectra are $\varepsilon$-close in Hausdorff distance to those of prescribed infinite covering graphs; they also gave a probabilistically strongly explicit variant.

A broader literature developed locally computable expanders by graph products and related transformations. Famously, in~\cite{RVW02},  Reingold--Vadhan--Wigderson introduced the \emph{zig-zag product}. Sooner after, Alon--Lubotzky--Wigderson \cite{ALW01} related this construction  to \emph{semidirect products}, Rozenman--Shalev--Wigderson \cite{RSW06} constructed expanding Cayley graphs on \emph{iterated wreath products}, and Rozenman--Vadhan \cite{RV05} developed \emph{derandomized squaring}. Cooper--Dotterrer--Prassidis \cite{CDP06} study covering properties of generalized zig-zag products. Ben-Aroya and Ta-Shma \cite{BATS11} obtained fully explicit almost-Ramanujan graphs with normalized nontrivial eigenvalues at most $d^{-1/2+o_d(1)}$ as $d\to\infty$, and, most recently, Jeronimo--Mittal--Roy--Wigderson \cite{JMRW25} give an operator-amplification transformation from arbitrary expanders to almost-Ramanujan expanders.

\smallskip\noindent\emph{Earlier explicit towers and automata.}
In~\cite[$\S$5]{Mak15}, Makelov constructed fully explicit expander towers using compatible graph products, including an almost-Ramanujan version based on \cite{BATS11}; his thesis also discusses binary refinements and existential Ramanujan $2$-lift towers.

Malyshev--Pak \cite{MP17} studied automaton-based lift towers and conjecture expansion for their level graphs. In \cite{BGV26},  Bondarenko--Grigorchuk--Vdovina constructed strongly explicit $(q+1)$-regular Ramanujan graphs from fixed Mealy automata for every odd prime power $q$, with $q$-fold covering levels.

Most recently, Kharlampovich--Vdovina identified the Bellaterra graphs with a cubic arithmetic Ramanujan tower in \cite{KV26a}. Their companion paper \cite{KV26b} proved expansion for the degree-six undirected Aleshin graphs and shows that they are not Ramanujan for $n\geq10$. These are distinct resolutions of the expansion conjectures of Malyshev--Pak. The present construction retains the field, quaternion algebra, and congruence conditions at~$2$ from the Bellaterra model and varies its adjacency prime.


\subsection{Algorithmic motivation}

The uses of expanders in theoretical computer science are surveyed by Hoory--Linial--Wigderson \cite{HLW06} and Vadhan \cite{Vad12}. For the present discussion, it is useful to separate expansion, sharp Ramanujan spectral bounds, and efficient local access to a compatible family of graphs.

\smallskip\noindent\emph{Randomness-efficient sampling.}
Expander walks replace independent samples by a short correlated walk; see Impagliazzo--Zuckerman \cite{IZ89} and Gillman's concentration theorem \cite{Gil98}. When the vertices are $n$-bit strings, strong explicitness lets the walk be generated without constructing a graph of size $2^n$. A compatible tower also allows a fine-level walk to project exactly to the corresponding coarse-level walk; this is proved in \cref{sec:consequences}.

\smallskip\noindent\emph{Codes, protographs, and communication networks.}  Famously,
Sipser--Spielman \cite{SS96} used expanders to construct asymptotically good codes with efficient decoding. Repeated lifts are particularly natural for protograph constructions: Ma--Yang \cite{MY07} study LDPC codes built by $2$-lifts, including stopping sets and implementation constraints. In \cite{VSDS16}, Valadarsky--Shahaf--Dinitz--Schapira used graph lifts in the Xpander datacenter architecture. These examples motivate locally accessible covering hierarchies. 

\smallskip\noindent\emph{Applications sensitive to sharp spectral information.}
In a surprizing development, Kahale \cite{Kah95} proved refined small-induced-subgraph estimates for Ramanujan graphs. Alon--Capalbo \cite{AC02} used Ramanujan graphs in constructions of unique-neighbor expanders, and Asherov--Dinur \cite{AD23} developed this direction for unbalanced bipartite graphs. Sharp spectra also control nonbacktracking mixing of random walks, see \cite{ABLS07} by Alon--Benjamini--Lubetzky--Sodin.  Curiously, in \cite{LP16},  Lubetzky--Peres established cutoff phenomenon for random walks on Ramanujan graphs and on graphs with vanishing spectral excess.

\smallskip\noindent\emph{Other algorithmic uses.}
Non-bipartite Ramanujan graphs, after scaling their edge weights, give standard spectral sparsifiers of the complete graph; see Batson--Spielman--Srivastava \cite{BSS12}. Expanders also enter Dinur's gap-amplification proof of the PCP theorem \cite{Din07}. Arithmetic graph constructions can also provide additional structure for quantum gate synthesis, as in Parzanchevski--Sarnak \cite{PS18}. These examples illustrate the range of settings in which explicit access to highly expanding graphs can be useful. 


\subsection{The construction and outline of the proof}\label{sec:outline}

Fix a rational prime $p\equiv1\pmod {17}$. The arithmetic takes place in the quadratic field $F=\Q(\sqrt{17})$ and its \emph{ring of integers}
\[
 \OF=\Z[\omega]=\{a+b\omega:a,b\in\Z\},\qquad
 \omega=\frac{1+\sqrt{17}}2\,.
\]
The ring of integers consists of the elements of $F$ satisfying a monic polynomial with integer coefficients. We distinguish \emph{rational primes}, which are integers (the usual prime numbers), from \emph{prime ideals} of this ring. An \emph{ideal} is an additive subgroup closed under multiplication by elements of the ring.  A nonzero proper ideal $R\subset\OF$ is prime (i.e., a prime ideal) precisely when $\OF/R$ is a field. Saying that $R$ is \emph{above} a rational prime $\ell$ means $R\cap\Z=\ell\Z$, or equivalently that $R$ contains the ideal $\ell\OF$.

In the construction, two prime ideals have different roles. We fix $Q=(2,\omega)$ once and for all, while $P=(p,\omega-s_p)$ depends on $p$, where $s_p$ is a chosen root of $X^2-X-4$ modulo $p$. Parentheses here denote an ideal generated by the enclosed elements: $(a,b)=a\OF+b\OF$. The quotients are
\[
 \OF/Q\simeq\F_2,\qquad \OF/P\simeq\F_p.
\]
Here $\F_\ell=\Z/\ell\Z$ for a rational prime $\ell$. The size $\IdealNorm{R}=|\OF/R|$ is the \emph{norm} of a prime ideal. Arithmetic at $Q$ will supply binary precision and the covering maps; arithmetic at $P$ will supply $\IdealNorm{P}+1=p+1$ neighbors. Section~\ref{sec:arithmetic} verifies these assertions and gives the definitions used in the outline.

\smallskip\noindent\emph{Step 1: Identify the binary vertex sets.}
We work in a four-dimensional algebra $D$ over $F$, with explicitly specified quaternion multiplication. Inside it, an \emph{order} is an $\OF$-subring containing $1$ that is also a full $\OF$-lattice: a finitely generated $\OF$-module spanning $D$ over $F$. In our construction it has an eight-element integer basis and plays the role of a lattice closed under multiplication. The orders $\EE_n$ are defined by increasingly restrictive matrix congruences at $Q$, with the conditions at every other prime ideal unchanged. For an order $\EE_n$, the set $\Cl(\EE_n)$ consists of suitable right ideals, identified when one is obtained from another by multiplication on the left by a nonzero quaternion. These are the arithmetic objects that will serve as vertices, not the prime ideals $P$ and $Q$ themselves.

The base order $\EE_0$ has \emph{class number one}: its ideal-class set has one element. Its invertible elements are all scalars from $F$. These two facts remove the global equivalences and stabilizers, leaving a concrete set of local matrix cosets. The resulting bijections are
\[
 \Cl(\EE_n)\simeq K/K_n\simeq\Z/2^n\Z,
\]
where $K$ is a fixed group of matrices modulo scalars and $K_n$ is the subgroup fixing $0$ modulo $2^n$. Their definitions and the bijections are given in Section~\ref{sec:arithmetic}. They respect reduction from level $n+1$ to level $n$. This is what proves the exact vertex count and the two-to-one maps; the fact $\IdealNorm{Q}=2$ alone would not prove either assertion.

\smallskip\noindent\emph{Step 2: Choose the $p+1$ neighbors once and for all.}
At $P$, the arithmetic neighbor choices correspond to lines (one-dimensional subspaces) in~$\F_p^2$.
There are $(p^2-1)/(p-1)=p+1$ of them. They specify right ideals of $\EE_0$, and class number one says that each can be written as $x_j\EE_0$ for a quaternion $x_j$. Thus we obtain $p+1$ quaternions, not $p+1$ new vertices.

At $Q$, quaternion multiplication has a fixed $2\times2$ matrix representation $\splitmap$ over $\Q_2$. The selected generators have matrices with entries in the \emph{$2$-adic integers} $\Z_2$, which are compatible residues modulo $2,4,8,\ldots$. Quaternion conjugation fixes the scalar coordinate and negates the other three coordinates. Applying $\splitmap$ to the conjugate $\bar x_j$ gives
\[
 g_j=\splitmap(\bar x_j)=
 \begin{pmatrix}a_j&b_j\\c_j&e_j\end{pmatrix}\,,\quad
 f_j(t)=\frac{a_jt+b_j/2}{2c_jt+e_j}\,,\quad \text{where} \quad 1\le j \le p+1.
\]
The entries $a_j,e_j$ and the determinant are odd, while $b_j,c_j$ are even. Hence the denominator is invertible modulo every power of $2$, and the \emph{neighbor with label} $j\in [p+1]$ is defined as $f_j(t)\bmod 2^n$. The same formulas work for all~$n$. Conjugation pairs the labels so that for the reverse-edge-label $j^*$ (see Section~\ref{s:notation}), we have $f_{j^*}=f_j^{-1}$, making the graph undirected. Section~\ref{sec:generators} proves that summing these neighbor actions is exactly the \emph{Hecke operator} $T_P$, the arithmetic operator that sums over the $p+1$ ideal-neighbor choices.

\smallskip\noindent\emph{Step 3: Obtain the exact spectral bound.}
The Hecke interpretation lets us apply a theorem from the theory of modular forms. The relevant vector space is simply the space of complex-valued functions on $\Cl(\EE_n)$, also called the \emph{Brandt module}. The arithmetic class-number conditions imply that its only exceptional, one-dimensional Hecke component consists of constant functions. The \emph{Jacquet--Langlands correspondence} transfers every other Hecke eigenvalue to Hilbert cusp forms over $F$, and the \emph{Ramanujan--Petersson theorem} gives
\[
 |\lambda|\leq2\sqrt{\IdealNorm{P}}=2\sqrt p\,.
\]
Section~\ref{sec:spectrum} states precisely the modular-form input and explains its terminology. No modular forms or eigenvalues need to be computed by the graph algorithm; this theorem certifies the spectrum of the already specified neighbor maps.

\smallskip\noindent\emph{Step 4: Prove that consecutive levels are genuine $2$-covers.}
For a nonzero $2$-adic number $z$, write $v_2(z)$ for the exponent of $2$ dividing it. Direct subtraction in the fractional-linear formula gives
\[
 v_2(f_j(t)-f_j(u))=v_2(t-u)\qquad(t\ne u).
\]
Thus each neighbor map preserves agreement of binary digits, and also the first digit at which two inputs differ. It maps the two children of a residue modulo $2^n$ bijectively to the two children of its neighbor. Reduction therefore gives a label-preserving two-sheeted cover $\Gamma_{n+1,p}\to\Gamma_{n,p}$. If $y=f_j(t)\bmod2^n$, the next binary digit of $f_j(t)\bmod2^{n+1}$ determines whether the lifted edge is straight or crossed. Section~\ref{sec:lifts} gives the formula and checks consistency with the reverse-edge label.

\smallskip\noindent\emph{Step 5: Evaluate neighbors without constructing the graph.}
Each quaternion $x_j$ is stored by eight integers in a fixed basis. A finite search, depending on $p$ but not on $n$, finds these coordinates. A query computes only the selected matrix $g_j$ to the needed binary precision and evaluates its fractional-linear formula. The algebraic constants are computed by \emph{Hensel lifting}, which extends a simple polynomial root modulo $2^k$ to a root modulo $2^{k+1}$ one bit at a time. Modular inversion handles the odd denominator. Section~\ref{sec:complexity} gives an $O_p(n^3)$ bit-operation bound, including the extra output bit needed for a lift sign. No adjacency matrix, previous graph level, or graph-sized signing table is stored.

The proof thus separates the choice of neighbors, the spectral estimate, and the binary covering structure. The prime ideal $P$ controls the first, modular-form theory supplies the second, and the fixed prime ideal $Q$ organizes the third. The local algorithm evaluates the resulting formulas rather than searching for a signing at each level.

\smallskip

\smallskip\noindent\emph{Number-theoretic tools.}
We verify the explicit quaternion order, local splitting, unit calculation, and vertex coordinates from Kharlampovich--Vdovina \cite{KV26a} for completeness. The three principal external inputs are the \emph{definite Eichler mass formula} \cite{KiV10}, the \emph{Jacquet--Langlands correspondence} \cite{DV13}, and \emph{Blasius's Ramanujan--Petersson theorem} for Hilbert modular forms \cite{Bla06}. We also use standard ideal and completion theory, quadratic reciprocity, Minkowski's bound, and special-value formulas for zeta and Dirichlet $L$-functions; the relevant sources are given where these facts are used. Dirichlet's theorem supplies the infinitude of admissible degrees \cite{Dav00}. On the other hand, the local-to-global lattice step, the neighbor formulas, the binary-covering rule, and the local-query bound are established below.
\smallskip

\subsection{Comparison with LPS}\label{subsec:lps-comparison}
The spectral argument follows the same broad strategy as LPS, with a different arithmetic setting.  As in \cite{LPS88}, an arithmetic graph operator is related to modular forms, and a Ramanujan--Petersson bound gives the optimal graph-eigenvalue estimate. LPS use Hamilton quaternions over~$\Q$; suitable norm-$p$ elements give the $p+1$ generators. Their original proof relates graph-walk counts to representation numbers of powers $p^k$ by a quaternary quadratic form. The associated theta series is a classical modular form of weight $2$, and the cuspidal coefficient bound yields the spectral estimate \cite[\S4]{LPS88}. In fact, LPS explicitly credit Eichler and Igusa for the weight-$2$ input \cite[\S1]{LPS88}.

In our approach, the adjacency operator is the Hecke operator $T_{\adjprime}$ on a quaternionic ideal-class module over $F=\Q(\sqrt{17})$. The Jacquet--Langlands correspondence transfers its nonconstant part to Hilbert cusp forms of parallel weight $(2,2)$ \cite[Thm~3.9]{DV13}. Blasius's Ramanujan--Petersson theorem \cite[Thm~1]{Bla06} gives
\[
 \lambda=\beta+\gamma,\qquad |\beta|=|\gamma|=\sqrt p,
 \qquad\text{hence}\qquad |\lambda|\leq2\sqrt p.
\]
The binary covering structure is a separate congruence construction at the fixed prime ideal $Q$ above $2$, rather than a consequence of this spectral theorem.

The separate arithmetic feature in our approach is the dyadic level structure: class number one identifies the vertices with $\Z/2^n\Z$, while fixed algebraic transformations give the covers and locally computable signings. Those properties require the explicit arithmetic and algorithmic arguments above; they do not follow from the Ramanujan--Petersson theorem alone.

\subsection{Paper structure}
Section~\ref{s:notation} fixes the graph conventions and the computational model. Section~\ref{sec:arithmetic} develops the arithmetic needed to identify the vertices of $\Gamma_{n,p}$, and Section~\ref{sec:generators} constructs the neighbors. Section~\ref{sec:spectrum} proves the Ramanujan bound. The covering property is established separately in Section~\ref{sec:lifts}, and Section~\ref{sec:complexity} completes the proof by establishing strong explicitness. Section~\ref{sec:consequences} records consequences of the binary organization.

\bigskip\section{Notation}\label{s:notation}

We first specify the graph conventions and the local-query model. Section~\ref{sec:arithmetic} separately defines the arithmetic notation, including prime ideals, completions, quaternion orders, and ideal classes. In particular, the capital letters $P,Q$ below will denote ideals, not integers or graph vertices.

We use undirected labeled multigraphs $G=(V,E)$ specified by permutations. Let $V$ be the vertex set and $[d]:=\{1,\ldots,d\}$ the set of \emph{edge labels}. For each $j\in[d]$, a \emph{permutation} (bijection) $f_j:V\to V$ specifies the neighbor reached from $t\in V$ along the edge incidence $(t,j)$.
We use $j^*\in [d]$ to denote the label for traversing edge $(t,j)$ in the opposite direction, and refer to
$j^*$ as the \emph{reverse-edge label}. An involution $j\mapsto j^*$ satisfies $f_{j^*}=f_j^{-1}$.
Following \cite[\S2.1]{RVW02}, we use the \emph{rotation-map notation}:
\[
 \operatorname{Rot}_G(t,j):=(f_j(t),j^*),
 \qquad \operatorname{Rot}_G^2=\mathrm{id}.
\]
Here the reverse-edge label depends only on $j$, not on the vertex~$t$. Labels distinguish incidences, not necessarily distinct neighbors; multiplicities are retained. This is not a proper edge coloring. The corresponding \emph{adjacency operator}~$A=A_G$ is
\[
 [Ah](t):=\sum_{j=1}^d h(f_j(t)), \quad \text{where} \quad t\in V, \ \ h: V\to \C.
\]
A fixed point of an edge label satisfying $j=j^*$ contributes one to the diagonal and is a \emph{half-loop}. A fixed point of a pair of distinct reverse-edge labels contributes two, as an ordinary loop does. A \emph{label-preserving covering} from a graph with vertex set $\widetilde V$ and neighbor permutations $\widetilde f_j$ to one with vertex set $V$ and neighbor permutations $f_j$ is a surjective map $\rho:\widetilde V\to V$ such that
\[
 \rho(\widetilde f_j(v))=f_j(\rho(v))
 \qquad(v\in\widetilde V,\ j\in[d]),
\]
with the same reverse-edge involution in both graphs. The induced map of incidences $(v,j)\mapsto(\rho(v),j)$ is then a bijection on the incidences at each vertex and commutes with reversal. The cover is \emph{two-sheeted} if every fiber $\rho^{-1}(t)$ has two elements. We use this rotation-map convention throughout; it permits an ordinary edge upstairs to project to a half-loop downstairs.

For fixed degree $d$, a family $\{G_n=(V_n,E_n)\}$ is \emph{strongly explicit} if there is a deterministic algorithm
with the following properties.  Given $n\ge 1$, $t\in V_n$, and an \emph{edge label} \, $j\in[d]$, the algorithm
returns the corresponding neighbor $t'\in V_n$ in $\operatorname{poly}(\log |V_n|)$ bit operations.  Additionally, the
reverse-edge label $j^\ast$ is included in the fixed finite description for every $j\in[d]$.
That description may depend on the fixed family (in our case on $p$), but not on $n$. Graph-sized preprocessing, level-dependent advice, and a separately stored signing table at each level are not allowed.

\bigskip\section{The fixed arithmetic model}\label{sec:arithmetic}

We first define the arithmetic objects and explain which ones become graph vertices
and which ones define edges. The prime ideals $Q$ and $P$ encode two different congruence
calculations: binary precision at $Q$ and the $p+1$ neighbor choices at $P$. We then construct the quaternion orders and identify their ideal classes with the residues modulo $2^n$.  For background on number fields and quaternion ideals, see e.g.\ \cite{KiV10,Neu99}.

\subsection{The field, ideals, and norms}\label{subsec:ideals}
The quadratic field is $F=\Q(\sqrt{17})=\{a+b\sqrt{17}:a,b\in\Q\}$. Put
\begin{equation}\label{eq:field}
 r=\sqrt{17},\qquad \omega=\frac{1+r}{2},\qquad
 \OF=\Z[\omega],\qquad
 \eta=\frac{5+r}{2}=\omega+2,\qquad
 \rho=\frac{5-r}{2}=3-\omega.
\end{equation}
The ring $\OF$ is the ring of integers of $F$: its elements are exactly those satisfying a monic polynomial with integer coefficients. In particular, every element of $\OF$ has a unique expression $a+b\omega$ with $a,b\in\Z$, and
\[
 \omega^2-\omega-4=0.
\]
The two real embeddings $\tau_+,\tau_-:F\to\mathbb R$ send $r$ to $\sqrt{17}$ and $-\sqrt{17}$, respectively. The \emph{field norm} and \emph{field trace} of $a\in F$ are
\[
 \FieldNorm{a}=\tau_+(a)\,\tau_-(a),\qquad
 \Tr_{F/\Q}(a)=\tau_+(a)+\tau_-(a).
\]
An element is \emph{totally positive} when both embeddings are positive. A \emph{unit} of a ring is an element whose multiplicative inverse also belongs to that ring; the superscript $\times$ denotes its unit group.

An \emph{ideal} $I\subseteq\OF$ is an additive subgroup such that $aI\subseteq I$ for all $a\in\OF$. The notation $(a_1,\ldots,a_k)$ denotes $a_1\OF+\cdots+a_k\OF$; in particular, $(a)=a\OF$ is a \emph{principal ideal}. The product $IJ$ is the ideal generated by all products $xy$ with $x\in I$ and $y\in J$, and $I^m$ is its $m$-fold product. Congruence $a\equiv b\pmod I$ means $a-b\in I$.

Every nonzero ideal has finite index. We write its \emph{ideal norm} with parentheses, $\IdealNorm{I}:=|\OF/I|$. For a nonzero principal ideal this satisfies $\IdealNorm{a\OF}=|\FieldNorm{a}|$. We use $\IdealNorm{I}$ only for ideals, reserving $\FieldNorm{a}$ for the field norm of an element. A nonzero proper ideal $R$ is \emph{prime} if $ab\in R$ implies $a\in R$ or $b\in R$; here this is equivalent to $\OF/R$ being a finite field, called its \emph{residue field}. We say $R$ \emph{lies above} a rational prime $\ell$ if $R\cap\Z=\ell\Z$, equivalently if $(\ell)\subseteq R$. Factorization of ideals also expresses this as $R\mid(\ell)$. Nonzero ideals of $\OF$ factor uniquely into prime ideals; this is the factorization used below. Thus $p,2$ are integers, whereas $P,Q,Q'$ below are subsets of $\OF$, namely prime ideals.

\subsection{The prime ideals used for levels and adjacency}\label{subsec:two-primes}
The polynomial $X^2-X-4$ reduces modulo $2$ to $X(X-1)$. Substituting either root defines a surjective ring homomorphism $\OF\to\F_2$, so its kernel is a prime ideal. The two kernels are
\[
 \levelprime=(2,\omega),\qquad
 \otherlevelprime=(2,\omega-1),\qquad
 \OF/\levelprime\simeq\OF/\otherlevelprime\simeq\F_2.
\]
For example, $a+b\omega\equiv a\pmod Q$ but $a+b\omega\equiv a+b\pmod{Q'}$. These ideals are distinct. The identities $\eta\rho=2$ and $\FieldNorm{\eta}=\FieldNorm{\rho}=2$ give
\begin{equation}\label{eq:dyadic-factorization}
 \levelprime=(\eta),\qquad \otherlevelprime=(\rho),\qquad
 (2)=\levelprime\otherlevelprime,
 \qquad \IdealNorm{\levelprime}=\IdealNorm{\otherlevelprime}=2.
\end{equation}
For instance, $\eta\in Q$ and both $(\eta)$ and $Q$ have index $2$, proving $(\eta)=Q$. This is what it means for $2$ to \emph{split} in $F$. We fix $Q=(\eta)$ for the entire construction; we never switch between the two prime ideals as $n$ or $p$ changes. The adjective \emph{dyadic} means associated with a prime ideal above $2$.

Now let $p\equiv1\pmod {17}$ be a rational prime. Quadratic reciprocity implies that $17$ is a nonzero square modulo $p$, so $X^2-X-4$ has two distinct roots in $\F_p$. Let $s_p$ be the smaller root in $\{0,\ldots,p-1\}$ and define
\begin{equation}\label{eq:activeprime}
 \adjprime=\adjprime_p=(p,\omega-s_p),\qquad \IdealNorm{\adjprime}=p.
\end{equation}
Indeed, the map $a+b\omega\mapsto a+bs_p\pmod p$ has kernel $P$, so $\OF/P\simeq\F_p$. The other root gives a different prime ideal $P'$, and $(p)=PP'$. Choosing the smaller root fixes $P$ deterministically; scanning residues computes it in $O(p\operatorname{poly}(\log p))$ bit operations. Since $p\ne2$, we have $P\ne Q$.

The division of labor is now explicit: $Q$ will control successive binary congruences, while the residue field $\OF/P\simeq\F_p$ gives the vector space $\F_p^2$, whose $p+1$ one-dimensional subspaces index the edge labels. Neither $P$ nor $Q$ is a graph vertex. The graph vertices will be classes of right ideals in quaternion orders, introduced below.

\subsection{Completions and binary precision}\label{subsec:completion}
For each nonzero prime ideal $R$, write $v_R(a)$ for the exponent of $R$ in the prime-ideal factorization of the principal fractional ideal $(a)$, with $v_R(0)=+\infty$. A fractional ideal is a nonzero finitely generated $\OF$-submodule of $F$ that becomes an integral ideal after multiplication by a nonzero element of $\OF$; its factorization may have negative exponents. The \emph{completion} $F_R$ is the field obtained by completing $F$ for the distance
\[
 |a-b|_R=(\IdealNorm{R})^{-v_R(a-b)}.
\]
Thus closeness means agreement modulo a large power of $R$, not closeness in the real line. Its ring of integers, denoted $\mathcal O_{F,R}$ in this paper, consists of the elements with nonnegative valuation. Throughout the paper, this subscript denotes the \emph{completed} ring. An element is a \emph{local unit at $R$} if its valuation is zero, equivalently if it is invertible in $\mathcal O_{F,R}$. A generator of the ideal of elements with positive valuation is called a \emph{uniformizer}. A finite place of $F$ means one of these completions; the real places are the two real embeddings.

In the case $R=Q$, the completion is particularly concrete. Let $\Z_2$ be the ring of compatible sequences of residues modulo $2^k$, and let $\Q_2$ be its fraction field. The even root of $X^2-X-4$ modulo $2$ lifts uniquely to compatible roots modulo every $2^k$: its derivative $2X-1$ is odd. More generally, this unique extension of a root with odd derivative is the simple-root form of \emph{Hensel lifting}. Sending $\omega$ to that even root embeds $F$ in $\Q_2$ and realizes
\begin{equation}\label{eq:q-completion}
 F_Q\simeq\Q_2,\qquad \mathcal O_{F,Q}\simeq\Z_2,
 \qquad Q\mathcal O_{F,Q}=2\Z_2.
\end{equation}
The image of $F$ is dense since it contains $\Q$, and its valuation is the one defined by $Q$. In these coordinates $\rho$ is odd and $\eta=2/\rho$ has valuation one. Thus divisibility by $Q^k$ becomes divisibility by $2^k$. In particular,
\[
 \OF/Q^k\simeq\Z/2^k\Z.
\]
This explains binary precision. Identifying the \emph{graph} vertices with these residues still requires the ideal-class argument below.

For $z\in\Q_2^\times$, $v_2(z)$ is the usual exponent of $2$, extended from rational numbers; $|z|_2=2^{-v_2(z)}$. In $\Z_2$, an element is \emph{odd} precisely when it is a unit. Computing a $2$-adic integer to precision $2^k$ means computing its residue modulo $2^k$. We never store an arbitrary infinite digit sequence: every constant used below is specified by an algebraic equation together with a chosen root modulo $2$.

\subsection{Quaternions, orders, and their matrix representation}\label{subsec:quaternions}
The quaternion algebra is the four-dimensional $F$-vector space with basis $1,\qi,\qj,\qk$ and multiplication determined by
\begin{equation}\label{eq:quaternion}
 D=(-1,-3)_F=F\langle\qi,\qj\rangle,
 \qquad \qi^2=-1,\quad \qj^2=-3,\quad \qi\qj=-\qj\qi=\qk.
\end{equation}
Bold letters denote quaternion units, not edge labels. For $x=x_0+x_1\qi+x_2\qj+x_3\qk$, quaternion conjugation fixes $x_0$ and negates the other three coordinates: $\bar x:=x_0-x_1\qi-x_2\qj-x_3\qk$. Its \emph{reduced trace} is $\trd(x)=x+\bar x=2x_0$, and its \emph{reduced norm} is
\begin{equation}\label{eq:norm}
 \nrd(x)=x\bar x=x_0^2+x_1^2+3x_2^2+3x_3^2\in F.
\end{equation}
This norm is multiplicative, and $x^{-1}=\bar x/\nrd(x)$ for $x\ne0$. It is distinct from both the rational-valued field norm $a\mapsto\FieldNorm{a}$ and the integer-valued ideal norm $I\mapsto\IdealNorm{I}$. At each real embedding, the displayed quadratic form is positive definite. This is the meaning of \emph{totally definite} for $D$.

An $\OF$-\emph{lattice} in $D$ is a finitely generated $\OF$-submodule that spans $D$ over $F$. An \emph{order} is such a lattice that is also a subring containing $1$; a \emph{maximal order} is one contained in no larger order. In our explicit bases an order has rank four over $\OF$, hence rank eight over $\Z$. For a prime ideal $R$, the subscript $R$ means extension of coefficients to the completed ring or field:
\[
 D_R=D\otimes_F F_R,\qquad
 \mathcal A_R=\mathcal A\otimes_{\OF}\mathcal O_{F,R}
\]
for an order or lattice $\mathcal A$. Equivalently, take its $R$-adic closure inside $D_R$. The algebra is \emph{split at $R$} if $D_R\simeq M_2(F_R)$, where $M_2$ denotes all $2\times2$ matrices; otherwise the local quaternion algebra is a division algebra and is called \emph{ramified}. Splitting a quaternion algebra is a matrix-representation property, different from splitting a rational prime into prime ideals.

The following explicit order and local splitting are taken from the arithmetic model of Kharlampovich--Vdovina \cite[\S4]{KV26a}.  Since these facts are used repeatedly below, we include their verification for completeness.  Define
\begin{equation}\label{eq:order}
 u=\frac{1+\qj}{2},\qquad
 v=\frac{\qi-r}{2}+\frac{\qk+r\qj}{6},\qquad
 \OO=\OF\langle1,\qi,u,v\rangle.
\end{equation}
Here the brackets mean the $\OF$-span of the four displayed elements.

\begin{proposition}[{\rm Explicit order and the splitting at $Q$}]\label[proposition]{prop:explicit-order}
The lattice $\OO$ in \cref{eq:order} is a maximal $\OF$-order in $D$.  Its reduced discriminant is the unit ideal, so $D$ is split at every finite prime ideal.  At the fixed prime $Q$ above $2$, the map
\begin{equation}\label{eq:splitting}
 \splitmap(\qi)=\begin{pmatrix}1&1+z\\z&-1\end{pmatrix},\qquad
 \splitmap(\qj)=\begin{pmatrix}-1&-2\\2&1\end{pmatrix},
\end{equation}
where $z\in\Z_2$ is the even root of $z^2+z+2=0$, extends to an isomorphism
\[
 \splitmap:D_Q\xrightarrow{\sim}M_2(\Q_2),\qquad
 \splitmap(\OO_Q)=M_2(\Z_2).
\]
Under this isomorphism reduced norm and reduced trace become determinant and matrix trace.
\end{proposition}

\begin{proof}
We follow the finite calculations in \cite[\S4.2 and \S4.4]{KV26a}.  First, direct multiplication in $D$, using $r^2=17$, gives
\[
\begin{array}{lll}
 \qi^2=-1, & u^2=u-1, & v^2=-6-rv,\\
 \qi u=2r-\qi-ru+3v, & u\qi=-2r+2\qi+ru-3v,\\
 \qi v=11-r\qi-6u+rv, & v\qi=-12+6u-rv,\\
 uv=-r+\qi-v, & vu=r-\qi-ru+2v.
\end{array}
\]
Every coefficient belongs to $\OF$, so the $\OF$-lattice in \cref{eq:order} is closed under multiplication and is therefore an order.

Let $b=(1,\qi,u,v)$.  The matrix of the trace pairing $(x,y)\mapsto\trd(xy)$ in this basis is
\[
 \bigl(\trd(b_sb_t)\bigr)_{s,t}=
 \begin{pmatrix}
 2&0&1&-r\\
 0&-2&0&-1\\
 1&0&-1&-r\\
 -r&-1&-r&5
 \end{pmatrix},
 \qquad \det=-1.
\]
Thus its trace-discriminant ideal is the unit ideal.  To see directly why this implies maximality, suppose locally that $\OO_R$ were properly contained in another order $\OO'_R$.  Expressing a basis of $\OO_R$ in a basis of $\OO'_R$ multiplies the determinant of the trace pairing by the square of the change-of-basis determinant.  A proper finite-index inclusion would therefore make the trace determinant divisible by the maximal ideal of $\mathcal O_{F,R}$, contradicting the unit determinant above.  Hence $\OO$ is maximal.  The standard discriminant--ramification relation for quaternion orders then says that every finite ramified prime divides the reduced discriminant of a maximal order; see \cite[\S2]{KiV10}.  Since that discriminant is $(1)$, there are no finite ramified primes.  At the two real embeddings the norm form in \cref{eq:norm} is positive definite, so $D$ is ramified there and is totally definite.

It remains to verify the displayed splitting and its integral lattice.  In the embedding $F_Q\simeq\Q_2$ from \cref{eq:q-completion}, put $\alpha=\omega-2$, so
\[
 \alpha^2+3\alpha-2=0,\qquad \alpha\equiv0\pmod2,
 \qquad r=2\alpha+3.
\]
The even roots $\alpha$ and $z$ are uniquely determined by Hensel lifting.  If $I=\splitmap(\qi)$ and $J=\splitmap(\qj)$ are the two matrices in \cref{eq:splitting}, then the relation $z^2+z+2=0$ gives
\[
 I^2=-I_2,\qquad J^2=-3I_2,
 \qquad IJ=-JI.
\]
Thus the defining quaternion relations hold, so \cref{eq:splitting} extends to an $F_Q$-algebra homomorphism $D_Q\to M_2(\Q_2)$.  The basis calculation below will prove that this map is bijective and will identify its integral image.

The images of the remaining basis elements are
\[
 \splitmap(u)=\begin{pmatrix}0&-1\\1&1\end{pmatrix},
\]
\[
 \splitmap(v)=\frac13
 \begin{pmatrix}
 -4\alpha+z-4&-2\alpha+2z-2\\
 2\alpha+z+2&-2\alpha-z-5
 \end{pmatrix}.
\]
Since $\alpha$ and $z$ are even and $3$ is a $2$-adic unit, all entries lie in $\Z_2$.  Vectorize a $2\times2$ matrix by listing its four entries.  The determinant of the $4\times4$ matrix whose columns are the vectorizations of
$\splitmap(1),\splitmap(\qi),\splitmap(u),\splitmap(v)$ is
\[
 -z^2-z-1=1,
\]
again using $z^2+z+2=0$.  It is therefore a unit of $\Z_2$, so these four matrices form a $\Z_2$-basis of $M_2(\Z_2)$.  Consequently $\splitmap(\OO_Q)=M_2(\Z_2)$, and extending scalars proves that $\splitmap:D_Q\to M_2(\Q_2)$ is an isomorphism. Finally, under any splitting of a quaternion algebra, the reduced characteristic polynomial agrees with the ordinary matrix characteristic polynomial, so reduced trace and norm become matrix trace and determinant.
\end{proof}

This is the \emph{fixed splitting} used to turn quaternions into neighbor matrices.  All its entries have finite algebraic descriptions, and Hensel lifting computes them to any prescribed binary precision.

\subsection{Orders with binary level}\label{subsec:level-orders}
The graph level $n$ is implemented by a matrix divisibility condition, not by changing the field or quaternion algebra. For $n\geq0$, define $\EE_n$ by its completed orders at all prime ideals:
\begin{equation}\label{eq:En}
 (\EE_n)_R=\OO_R\quad(R\ne Q),\qquad
 (\EE_n)_Q=
 \begin{pmatrix}\Z_2&2^{n+1}\Z_2\\2\Z_2&\Z_2\end{pmatrix}.
\end{equation}
The matrix notation means that the upper-right entry is divisible by $2^{n+1}$, the lower-left entry by $2$, and the other entries are arbitrary $2$-adic integers. Globally, $\EE_n$ consists of the elements of $\OO$ satisfying this condition under $\splitmap$; the conditions at $R\ne Q$ are unchanged. The local-to-global lattice argument in \cref{lem:local-global} verifies that its completions are exactly the prescribed orders. Their intersection is closed under multiplication and contains $1$, so it is an order. Thus $\EE_{n+1}\subset\EE_n$.

An \emph{Eichler order} is the intersection of two maximal orders. At a split prime, its standard local form consists of integral matrices whose lower-left entry is divisible by a prescribed prime power; that prime power is its local \emph{level}. With $h=\diag(2^{n+1},1)$, where $\diag$ denotes a diagonal matrix, we have
\[
 h^{-1}(\EE_n)_Qh=
 \begin{pmatrix}\Z_2&\Z_2\\2^{n+2}\Z_2&\Z_2\end{pmatrix}.
\]
Consequently $\EE_n$ is an Eichler order of global level $Q^{n+2}$. The arithmetic level $Q^{n+2}$ is an ideal; it is not the graph level $n$ or the number $2^n$ of vertices. The shift by two is part of this fixed choice of base order. We use $n=0$ internally for the arithmetic base, although the graph family in the main theorem starts at $n=1$.

The base order has the global basis
\begin{equation}\label{eq:E0basis}
 \EE_0=\OF\langle1,\eta\qi,\eta u,v\rangle.
\end{equation}
To check it, reduce \cref{eq:splitting} modulo $2$. The upper/lower off-diagonal residues of $\qi,u,v$ are $(1,0),(1,1),(0,0)$. Thus the coefficients of $\qi$ and $u$ must belong to $Q$, whereas those of $1,v$ are unrestricted. At every $R\ne Q$, $\eta$ is a local unit. This proves \cref{eq:E0basis} from the completed conditions.

\subsection{Ideal classes, scalar units, and the vertex count}\label{subsec:class-definitions}
There are two different class-number notions in the proof. For the number field, the \emph{ordinary ideal class group} is the group of nonzero fractional ideals of $\OF$ modulo principal fractional ideals. The \emph{narrow class group} uses only principal ideals generated by totally positive elements. Narrow class number one therefore says that every nonzero ideal has a totally positive generator. This will let us choose $P=(\pi)$ with $\FieldNorm{\pi}=p$; it is not a statement about the number of graph vertices.

For a quaternion order $\mathcal A$, a \emph{right $\mathcal A$-ideal} $I$ is an $\OF$-lattice in $D$ satisfying $I\mathcal A\subseteq I$. It is \emph{integral} if $I\subseteq\mathcal A$, \emph{principal} if $I=x\mathcal A$ for some $x\in D^\times$, and \emph{locally principal} if $I_R=x_R\mathcal A_R$ for some $x_R\in D_R^\times$ at every finite prime ideal $R$. Here $D^\times$ and $D_R^\times$ denote the invertible elements. Two locally principal right ideals $I,J$ are equivalent if $J=xI$ for one $x\in D^\times$. Their finite set of equivalence classes is $\Cl(\mathcal A)$, and its size is the \emph{class number of the order}. In general this is a class \emph{set}, not a group. Its elements will be the graph vertices.

The \emph{left order} of $I$ is $\mathcal A_L(I)=\{x\in D:xI\subseteq I\}$. Its units are the invertible transformations preserving $I$. Scalar units are those lying in the central copy of $F$; they form $\OF^\times$. The quotient $\mathcal A_L(I)^\times/\OF^\times$ is the \emph{projective unit stabilizer} of the class, meaning that scalar transformations have been ignored. Finally, the reduced norm of an ideal is the fractional ideal
\[
 \nrd(I)=\bigl(\nrd(x):x\in I\bigr)\subset F.
\]
It is a fractional ideal of $\OF$, integral when $I$ is integral, whereas the reduced norm of an individual quaternion is an element of $F$. These conventions agree with \cite[\S1]{KiV10}.

\smallskip\noindent\emph{Passing between local and global lattices.}
We will use the following standard lattice fact, recording a proof to make the existence of the ideals used below explicit; see also \cite[\S1]{KiV10}.
\begin{lemma}[{\rm Local-to-global lattices}{}]\label[lemma]{lem:local-global}
For each finite prime ideal $R$, let $\Lambda_R\subset D_R$ be a full $\mathcal O_{F,R}$-lattice, and suppose that $\Lambda_R=\OO_R$ for all but finitely many $R$. There is a unique full $\OF$-lattice $\Lambda\subset D$ with these completions, namely
\[
 \Lambda=D\cap\prod_R\Lambda_R.
\]
If every $\Lambda_R$ is a right ideal of $\mathcal A_R$ for a fixed order $\mathcal A$, then $\Lambda$ is a right $\mathcal A$-ideal.
\end{lemma}
\begin{proof}
Choose a nonzero integer $a$ such that
$a\OO_R\subseteq\Lambda_R\subseteq a^{-1}\OO_R$ for every $R$. Such a choice is possible because only finitely many completed lattices differ from $\OO_R$. The finite $\OF$-module $a^{-1}\OO/a\OO$ is annihilated by $a^2$. The Chinese remainder theorem decomposes it into its primary components, which are precisely the nonzero quotients $a^{-1}\OO_R/a\OO_R$. In this decomposition, take the inverse image of the product of the submodules $\Lambda_R/a\OO_R$. This defines a lattice between $a\OO$ and $a^{-1}\OO$ with the prescribed completions. The same decomposition shows that membership is determined by these completed conditions, proving both the displayed intersection formula and uniqueness. If the completed lattices are right ideals, this intersection is stable under right multiplication by every element of $\mathcal A$.
\end{proof}

\begin{proposition}[{\rm Arithmetic base data}{}]\label[proposition]{prop:class}
The field $F$ has narrow class number one. Every Eichler order of level $Q^m$ with $m\geq2$, including the left orders of its locally principal right ideals, has only scalar units. Moreover
\begin{equation}\label{eq:classnumber}
 |\Cl(\EE_n)|=2^n,\qquad \EE_0^\times=\OF^\times.
\end{equation}
\end{proposition}
\begin{proof}
This is the calculation of \cite[Lemma~7.1]{KV26a}, rewritten in our notation; we include the proof for completeness.

First consider the class groups of the field.  Minkowski's bound for the real quadratic field of discriminant $17$ says that every ordinary ideal class contains an integral ideal of norm at most
$\sqrt{17}/2<2.1$.  Hence only norms $1$ and $2$ can occur.  The two prime ideals of norm $2$ are $(\eta)=Q$ and $(\rho)=Q'$, both principal, so the ordinary class number is one.

We next determine the signs of generators.  Put $\fundunit=4+r$; then $\FieldNorm{\fundunit}=-1$.  We record the elementary unit calculation from \cite[Lemma~7.1]{KV26a}.  Suppose a unit $w$ has first real embedding satisfying $1<w<\fundunit$.  Since $\FieldNorm{w}=\pm1$, its other real value is $w'=\pm w^{-1}$.  Write $w=(a+br)/2$ with $a,b\in\Z$ of the same parity.  Then
\[
 0<br=w-w'<2r,
\]
where one uses $\fundunit+\fundunit^{-1}=2r$ and
$\fundunit-\fundunit^{-1}=8<2r$.  Thus $b=1$.  But the equation
$\FieldNorm{w}=\pm1$ would then give $a^2=17\pm4$, i.e. $a^2=13$ or $21$, impossible.  Therefore there is no unit with first embedding strictly between $1$ and $\fundunit$.  Dividing any unit positive at the first embedding by a suitable power of $\fundunit$ places it in $[1,\fundunit)$, so every such unit is a power of $\fundunit$.  Since the second embedding of $\fundunit$ is negative, the totally positive units are exactly $\fundunit^{2k}$, hence are squares.  The sign patterns furnished by $-1$ and $\fundunit$ show that every principal ideal has a totally positive generator.  Since the ordinary class number is one, the narrow class number is also one.

Now let $\mathcal R$ be an Eichler order of level $Q^m$ with $m\ge2$, or the left order of a locally principal right ideal of such an order.  The latter has locally conjugate Eichler orders and hence the same local level.  Let $x\in\mathcal R^\times$.  Its reduced norm is a totally positive unit, and the preceding paragraph shows that it is a square of a scalar unit.  After multiplying $x$ by the inverse of that scalar unit, we may assume $\nrd(x)=1$.
If $x$ is noncentral, then at each real embedding the norm-one quaternion has reduced trace strictly between $-2$ and $2$.  Write
$t=\trd(x)=(a+br)/2\in\OF$, with $a,b\in\Z$ of the same parity.  The difference of the two real values of $t$ is $br$, so $|b|r<4$ and therefore $b=0$.  Hence $t\in\Z$ and $t\in\{-1,0,1\}$.

On the other hand, at $Q$ the order is locally conjugate to the standard Eichler order of level $2^m$, and $m\ge2$.  Modulo $4$, a determinant-one unit in that order is upper triangular with odd diagonal entries whose product is $1$.  Its two diagonal entries therefore have sum $2\pmod4$.  Thus $\trd(x)\equiv2\pmod4$, contradicting $t\in\{-1,0,1\}$.  Every unit is therefore central.  This proves both the scalar-unit assertion for $\mathcal R$ and, in particular, $\EE_0^\times=\OF^\times$.

It remains to count ideal classes. We apply the definite Eichler mass formula \cite[\S5, equations (5.1)--(5.2)]{KiV10}, together with the standard special-value identities specified below.  Because all projective unit stabilizers have just been proved trivial, the mass equals the ordinary number of ideal classes.  Let $\chi(a)=(a/17)$ be the quadratic character modulo $17$.  The nonzero quadratic residues modulo $17$ are
\[
 1,2,4,8,9,13,15,16,
\]
and direct summation gives
\[
 \sum_{a=1}^{16}\chi(a)a^2=136.
\]
Here $L(s,\chi)$ is the Dirichlet $L$-function, $\zeta_F(s)$ is the Dedekind zeta function of $F$, and $\zeta(s)$ is the Riemann zeta function. Initially, for $\operatorname{Re}s>1$, they are defined by
\[
 L(s,\chi)=\sum_{k\geq1}\frac{\chi(k)}{k^s},\qquad
 \zeta_F(s)=\sum_{0\ne I\subseteq\OF}\frac1{\IdealNorm{I}^{s}},
 \qquad \zeta(s)=\sum_{k\geq1}\frac1{k^s}.
\]
The sum defining $\zeta_F$ is over nonzero integral ideals. Their values at $s=-1$ refer to analytic continuation. The standard identity
\[
 L(-1,\chi)=-\frac{17}{2}\sum_{a=1}^{16}\chi(a)B_2(a/17),
 \qquad B_2(x)=x^2-x+\frac16,
\]
follows from the Hurwitz-zeta identities \cite[equations (25.11.14), (25.15.3)]{DLMF}. Since $\sum_a\chi(a)=\sum_a a\chi(a)=0$ for this even nonprincipal character, it gives
\[
 L(-1,\chi)=-\frac1{2\cdot17}\sum_{a=1}^{16}\chi(a)a^2=-4.
\]
The quadratic-field factorization $\zeta_F(s)=\zeta(s)L(s,\chi)$ and $\zeta(-1)=-1/12$ now yield
\[
 \zeta_F(-1)=\frac13.
\]
The finite discriminant of $D$ is $(1)$ by \cref{prop:explicit-order}, the ordinary class number of $F$ is one, $[F:\Q]=2$, and
$\IdealNorm{Q^m}=2^m$.  The mass formula therefore gives
\[
 |\Cl(\mathcal R)|
 =\frac12\cdot\frac13\cdot2^m
   \left(1+\frac1{\IdealNorm{Q}}\right)
 =\frac12\cdot\frac13\cdot2^m\cdot\frac32
 =2^{m-2}.
\]
Taking $m=n+2$ yields $|\Cl(\EE_n)|=2^n$.  This proves every assertion of the proposition.  Standard background on Minkowski's theorem, units, and zeta functions is reviewed in \cite{Neu99}.
\end{proof}

Class number one for $\EE_0$ means that every locally principal right ideal of $\EE_0$ is generated by one global element. This is why the neighbor ideals can be represented by finitely many quaternions. The separate assertion about scalar units ensures that a normalized generator is unique up to sign and that the class-set coordinates below have no further stabilizer identifications.

\subsection{Local coordinates for every vertex}\label{subsec:vertex-coordinates}
We now describe explicitly why the ideal classes form the desired binary tower. The notation used to assemble local information is as follows:
\[
 \widehat{\EE}_n=\prod_R(\EE_n)_R,\qquad
 \mathcal U_n=\prod_R(\EE_n)_R^\times.
\]
The product is over all nonzero prime ideals. The group $\widehat D^\times$ consists of tuples $(g_R)_R$ with $g_R\in D_R^\times$ and $g_R\in\OO_R^\times$ for all but finitely many $R$. Such a tuple is called a \emph{finite quaternionic idele}; the hat collects all finite places, not just $Q$. A global quaternion embeds diagonally, by using the same element in every completion. The double quotient $D^\times\backslash\widehat D^\times/\mathcal U_n$ identifies $g$ and $xgu$ for $x\in D^\times$ and $u\in\mathcal U_n$. The tuple $g$ represents the right ideal $I(g)=D\cap g\widehat{\EE}_n$. By \cref{lem:local-global}, its completion at $R$ is $g_R(\EE_n)_R$. Conversely, local generators of a locally principal right ideal form such a tuple, chosen to be units at all but finitely many places. Finally, $I(g')=xI(g)$ holds precisely when $g'=xgu$ for some $u\in\mathcal U_n$. This proves the identification with $\Cl(\EE_n)$ rather than merely assigning a formal double quotient to it.

For a commutative ring $S$ used here, $\GL_2(S)$ is the group of invertible $2\times2$ matrices and $\PGL_2(S)$ is its quotient by invertible scalar matrices. Ignoring those scalars does not change the fractional-linear action. Put
\[
 K=\ker\bigl(\PGL_2(\Z_2)\longrightarrow\PGL_2(\F_2)\bigr),
\]
where the map reduces matrix entries modulo $2$. Thus representatives of elements of $K$ have odd diagonal entries and even off-diagonal entries. Let $K_n$ be the image of $(\EE_n)_Q^\times$ in $\PGL_2(\Z_2)$. For $g=\left(\begin{smallmatrix}a&b\\c&e\end{smallmatrix}\right)\in K$, the usual fractional-linear action in the coordinate $x=2t$ becomes
\begin{equation}\label{eq:mobius}
 \mobius_g(t)=\frac{at+b/2}{2ct+e}.
\end{equation}
The denominator is odd, so this defines a map on $\Z_2$ and on every $\Z/2^n\Z$. The translations $\left(\begin{smallmatrix}1&2t\\0&1\end{smallmatrix}\right)$ send $0$ to $t$, proving transitivity. The condition that $\mobius_g(0)=0\pmod{2^n}$ is $b\in2^{n+1}\Z_2$, exactly the condition defining $K_n$. Thus $K_n$ is the stabilizer of $0$, and $K/K_n$ is the set of cosets $gK_n$, not necessarily a quotient group.

The class-set identification in the next lemma is the same local-to-global calculation used in \cite[Proposition~7.2]{KV26a}; we include the proof for completeness because it is what turns the arithmetic level into exactly $2^n$ graph vertices.

\begin{lemma}[{\rm All vertices, with compatible coordinates}{}]\label[lemma]{lem:vertices}
There are compatible bijections
\begin{equation}\label{eq:vertices}
 \Cl(\EE_n)
 =D^\times\backslash\widehat D^\times/\mathcal U_n
 \simeq\mathcal U_0/\mathcal U_n
 \simeq K/K_n
 \simeq\Z/2^n\Z.
\end{equation}
\end{lemma}
\begin{proof}
Class number one at level zero gives $\widehat D^\times=D^\times\mathcal U_0$. The remaining left equivalence on $\mathcal U_0/\mathcal U_n$ is by $D^\times\cap\mathcal U_0=\EE_0^\times=\OF^\times$. It is trivial because these elements are central and belong to every $\mathcal U_n$. Only the local factor at $Q$ changes with $n$, giving the coset set $K/K_n$.

The last bijection is $gK_n\mapsto\mobius_g(0)\bmod2^n$. An explicit idele representing the vertex $t$ has component $\left(\begin{smallmatrix}1&2t\\0&1\end{smallmatrix}\right)$ at $Q$ and identity component everywhere else. The map from level $n+1$ to level $n$ forgets the stronger subgroup condition; in these coordinates it is exactly reduction modulo~$2^n$.
\end{proof}

We can therefore take the graph vertices to be $t\in\Z/2^n\Z$. The two vertices above $t$ at the next level are $t$ and $t+2^n$. The bijection above establishes both the exact vertex count and its compatibility across levels.

This completes Step~1 of the proof outlined in Section~\ref{sec:outline}.

\bigskip\section{A fixed finite set of neighbor transformations}\label{sec:generators}

Each quaternion generator specifies a neighbor function, rather than a graph vertex. We construct the $p+1$ generators from arithmetic ideal neighbors and use their $2\times2$ matrix representations to obtain fractional-linear maps on $\Z_2$. Identifying the sum of these neighbor actions with a Hecke operator then connects the concrete graph construction to the spectral theorem in the next section.

For the prime ideal $\adjprime=\adjprime_p$ in \cref{eq:activeprime}, choose a totally positive generator $\pi$: this means $P=\pi\OF$ and $\tau_+(\pi),\tau_-(\pi)>0$. Then $\FieldNorm{\pi}=p$. Existence follows from narrow class number one, defined and proved in \cref{subsec:class-definitions}. A \emph{norm-$P$ ideal neighbor of $\EE_0$} is an integral locally principal right ideal $J\subset\EE_0$ with $\nrd(J)=P$. The elements used to generate these ideals are selected by the equation
\begin{equation}\label{eq:S}
 S_p=\{x\in\EE_0:\nrd(x)=\pi\}/\{\pm1\}.
\end{equation}
The quotient by $\{\pm1\}$ identifies $x$ and $-x$; we write $[x]$ for this sign class. This is different from an ideal class $[I]\in\Cl(\EE_n)$. Section~\ref{sec:complexity} gives an algorithm for computing $S_p$.

\begin{lemma}[{\rm Norm generators}{}]\label[lemma]{lem:normgenerators}
The set $S_p$ has exactly $p+1$ elements. Its image in $D^\times/F^\times$, where quaternions differing by a nonzero scalar are identified, is closed under taking inverses. Every representative is invertible in $(\EE_0)_Q$ and has scalar matrix reduction modulo $2$ under $\splitmap$.
\end{lemma}
\begin{proof}
Integral locally principal right $\EE_0$-ideals of reduced norm $\adjprime$ are indexed by the lines in $\F_p^2$. Indeed, at $\adjprime$ the order is a full matrix ring, and such a right ideal consists of matrices whose columns reduce into a prescribed line. Away from $\adjprime$ it is the full local order. By \cref{lem:local-global}, each such system of local right ideals gives a unique global right ideal, which is locally principal. There are $p+1$ choices, since $\F_p^2$ has $p^2-1$ nonzero vectors and each one-dimensional subspace has $p-1$ nonzero vectors. Here a line means a one-dimensional linear subspace, not an affine line.

Every one of these ideals is principal because $|\Cl(\EE_0)|=1$. Write it as $x\EE_0$. Its generator lies in $\EE_0$, and $\nrd(x)=\pi\varepsilon$ for a totally positive unit $\varepsilon$. Since $\varepsilon$ is a square, scalar rescaling makes $\nrd(x)=\pi$.

If two such normalized generators give the same right ideal, they differ on the right by a unit of $\EE_0$, hence a scalar unit $u$. Their equal reduced norms imply $u^2=1$, so they differ only by sign. This proves the count.

Conjugation preserves $\EE_0$, since $\bar x=\trd(x)-x$ and the reduced trace is integral. It also preserves the norm equation. As $\bar x=\pi x^{-1}$, it induces projective inversion on $S_p$. Finally, $\pi$ is a unit at $\levelprime$, so $x\in(\EE_0)_{\levelprime}^{\times}$, whose reduction is scalar modulo $2$.
\end{proof}

Choose representatives $x_1,\ldots,x_d$, with $d=p+1$, for these sign classes in a fixed order. To turn each quaternion into a neighbor rule, apply the matrix representation $\splitmap$ from \cref{eq:splitting} to its conjugate and put
\begin{equation}\label{eq:gj}
 g_j=\splitmap(\bar x_j)=\begin{pmatrix}a_j&b_j\\c_j&e_j\end{pmatrix},
 \qquad f_j(t)=\frac{a_jt+b_j/2}{2c_jt+e_j}.
\end{equation}
The matrices $g_j$ are in $\GL_2(\Z_2)$ and have odd determinant $\pi$. Their reduction modulo $2$ is scalar, so $a_j,e_j$ are odd and $b_j,c_j$ are even. In particular, $b_j/2$ is integral and $2c_jt+e_j$ is odd for every $t\in\Z_2$.

The formula for $f_j$ is the usual fractional-linear matrix action, written in the coordinate $2t$. Indeed,
\[
 g_j\begin{pmatrix}2t\\1\end{pmatrix}
 =\begin{pmatrix}2a_jt+b_j\\2c_jt+e_j\end{pmatrix}
 =(2c_jt+e_j)\begin{pmatrix}2f_j(t)\\1\end{pmatrix}.
\]
Dividing by the invertible bottom entry gives the affine coordinate $2f_j(t)$. Thus the neighbor of $t\in\Z/2^n\Z$ with edge label $j$ is
\[
 f_j(t)\bmod2^n
 =(a_jt+b_j/2)(2c_jt+e_j)^{-1}\pmod{2^n}.
\]
Here the inverse is taken modulo $2^n$. Replacing $x_j$ by $-x_j$ changes neither the fractional-linear map nor the neighbor, so the sign classes used above suffice. The parity conditions and the odd determinant make each map a permutation modulo every power of $2$.

The same quaternion $x_j$ therefore specifies the $j$th neighbor of every vertex at every level. Only the vertex $t$ and modulus $2^n$ change. The $p+1$ sign classes give $p+1$ edge labels, not necessarily $p+1$ distinct neighboring vertices; coincident images retain their multiplicities. Define
\begin{equation}\label{eq:adjacency}
 (A_{n,p}h)(t)=\sum_{j=1}^{d}h\bigl(f_j(t)\bmod2^n\bigr).
\end{equation}
Quaternion conjugation supplies the reverse-edge labels. The identity $\bar x_j=\pi x_j^{-1}$ shows that conjugation agrees with inversion up to a scalar, and scalars do not change a fractional-linear action. Hence the involution $j\mapsto j^*$ defined by $[x_{j^*}]=[\bar x_j]$ satisfies $f_{j^*}=f_j^{-1}$. Following edge label $j$ and then reverse-edge label $j^*$ returns to the original vertex. The resulting rotation map is
\begin{equation}\label{eq:rotation}
 \operatorname{Rot}_{n,p}(t,j)
 =\bigl(f_j(t)\bmod2^n,j^*\bigr).
\end{equation}
Thus \cref{eq:adjacency} is symmetric and has row sum $d$. The reverse-edge involution is independent of the level.

The norm condition is what distinguishes these neighbor rules from an arbitrary family of permutations. The $p+1$ quaternions represent the arithmetic ideal neighbors at $\adjprime$, so the sum of their neighbor actions is a Hecke operator. The next proposition makes this identification precise; it is the link between the explicit neighbor formulas and the Ramanujan bound proved in \cref{sec:spectrum}.

To define the Hecke operator independently of our coordinates, take a class $[I]\in\Cl(\EE_n)$. Its $P$-neighbors are the locally principal right ideals $J\subset I$ such that $J_R=I_R$ for $R\ne P$ and $\nrd(J)=P\nrd(I)$. At $P$ the order is a full matrix ring, so there are $p+1$ such ideals, indexed by the lines of $\F_p^2$. Define
\[
 (T_P h)([I])=\sum_{J\text{ a }P\text{-neighbor of }I}h([J]).
\]
Changing the representative from $I$ to $xI$ bijects these neighbors by $J\mapsto xJ$: the completed containments are preserved and $\nrd(xJ)=\nrd(x)\nrd(J)$. Thus $T_P$ is well-defined on ideal classes. The resulting classes need not be distinct; each ideal contributes once. The operator is \emph{unnormalized} because it is a sum, not an average: $T_P\one=(p+1)\one$. The name \emph{Hecke correspondence} refers to this neighbor relation, and its branches are the individual choices of $J$. This explains both the degree and the multiplicity convention in the next proposition; see \cite[\S4]{DV13}.

\begin{proposition}[{\rm Identification with the Hecke operator}{}]\label[proposition]{prop:hecke}
Under \cref{eq:vertices}, $A_{n,p}$ is the unnormalized Hecke operator $T_{\adjprime}$ on the \emph{entire} ideal-class set $\Cl(\EE_n)$.
\end{proposition}
\begin{proof}
Put $U=\GL_2(\mathcal O_{F,P})$. The local neighbor matrices form the double coset
\[
 U\diag(\pi,1)U=\coprod_{j=1}^{p+1}(x_j)_P U
 \quad\text{inside }\GL_2(F_P).
\]
Here $(x_j)_P$ denotes the matrix of $x_j$ in a splitting at $P$, and $\coprod$ denotes disjoint union. A double coset $UgU$ allows multiplication of $g$ by elements of $U$ on both sides; each $(x_j)_P U$ is one coset modulo multiplication on the right by $U$. Their distinctness is the distinctness of the norm-$P$ right ideals in \cref{lem:normgenerators}.

Let $t_j$ be the finite quaternionic idele (the tuple of local invertible elements defined in \cref{subsec:vertex-coordinates}) equal to $x_j$ at $P$ and to $1$ elsewhere. It factors as $t_j=x_jz_j$, where $(z_j)_{\adjprime}=1$ and $(z_j)_{\localprime}=x_j^{-1}$ at every other finite place. Since $x_j$ has unit norm away from $\adjprime$, $z_j\in\mathcal U_0$.

Represent a vertex by an idele $u$ whose components are the identity except possibly at $Q$, as in \cref{lem:vertices}. Because $\adjprime\ne\levelprime$, $u$ commutes with $t_j$. Its $j$th Hecke neighbor is therefore
\[
 [ut_j]=[t_ju]=[x_jz_ju]=[z_ju].
\]
In $K/K_n$ this is left multiplication by $\splitmap(x_j^{-1})$, which has the same projective action as $g_j=\splitmap(\bar x_j)$. Thus it is exactly $f_j$. This argument identifies each edge label, including coincident neighbors and diagonal contributions, not merely the total number of vertices.
\end{proof}

\begin{remark}
The induced local map uses $x_j^{-1}$, rather than $x_j$. Formula \cref{eq:gj} fixes the direction. In the cubic case the projective generators are involutions, so this distinction is invisible.
\end{remark}

This completes Step~2 of the proof outlined in Section~\ref{sec:outline}.

\bigskip\section{The exact spectral bound}\label{sec:spectrum}

The graph and its neighbor rules are already defined; this section proves their spectral bound. We explain the function space on which the Hecke operators act and the modular-form theorem applied to it. The arithmetic class-number conditions ensure that constants are the only exceptional component, so the bound controls every eigenvector orthogonal to them.

\smallskip\noindent\emph{The function space and the modular-form input.}
The \emph{Brandt module} in this proof is simply
\[
 \mathcal M_n=\{h:\Cl(\EE_n)\to\C\},
\]
with the Hecke operators defined by sums over ideal neighbors. The word ``module'' refers to this vector space together with those linear operators; its vectors are ordinary functions on the vertices. Write $\one$ for the constant function $1$, and $\one^\perp$ for the functions whose sum over the vertices is zero. The inner product is $\langle h_1,h_2\rangle=\sum_t h_1(t)\overline{h_2(t)}$.

For clarity, the analytic objects used by the spectral theorem are \emph{Hilbert cusp forms of parallel weight $(2,2)$}. They are holomorphic functions of two upper-half-plane variables, one for each real embedding of $F$. If $\mathcal H=\{z\in\C:\operatorname{Im}z>0\}$ and $\mathcal N$ is an ideal, the level group can be written
\[
 \Gamma_0(\mathcal N)=
 \left\{\begin{pmatrix}a&b\\c&e\end{pmatrix}\in
 \operatorname{SL}_2(\OF):c\in\mathcal N\right\}.
\]
Here $\operatorname{SL}_2$ means determinant-one matrices. For $\gamma$ in this group, let $a_\nu,b_\nu,c_\nu,e_\nu$ be its entries under $\tau_\nu$, for $\nu\in\{+,-\}$, and let $\gamma_\nu z=(a_\nu z+b_\nu)/(c_\nu z+e_\nu)$. The weight condition is
\[
 f(\gamma_+z_+,\gamma_-z_-)
 =(c_+z_++e_+)^2(c_-z_-+e_-)^2f(z_+,z_-).
\]
Thus ``parallel'' means the same weight $2$ in both variables. The cusp condition requires the constant term of the Fourier expansion to vanish at every cusp, that is, at each boundary class represented by a point of $F\cup\{\infty\}$.

The standard Hilbert modular description also allows integral matrices with totally positive unit determinant. In this field every totally positive unit is a square by \cref{prop:class}. If such a matrix has determinant $u^2$ with $u\in\OF^\times$, multiplying it by the scalar $u^{-1}$ gives a determinant-one matrix with the same fractional-linear action and the same level condition. For parallel weight $(2,2)$ the scalar change contributes $\FieldNorm{u}^{-2}=1$, so the displayed transformation law is unchanged. This explains the equivalent projective level-group description. We use the Hecke module with \emph{trivial central character}, meaning that scalar ideles act trivially; this central action is checked separately in the proof below. See \cite[\S\S2--3]{DV13}. The graph algorithm does not evaluate these analytic functions.

The \emph{Jacquet--Langlands correspondence} identifies the relevant quaternionic Hecke eigenvalues with those on Hilbert cusp forms, preserving their values at prime ideals away from the level. The \emph{Ramanujan--Petersson theorem} bounds these eigenvalues. For an eigenvalue $\lambda_P$ at a prime ideal $P$ of norm $p$, the two roots $\beta,\gamma$ of
\[
 T^2-\lambda_P T+p
\]
are called its arithmetic \emph{Satake parameters}. The theorem asserts $|\beta|=|\gamma|=\sqrt p$, so $|\lambda_P|\leq2\sqrt p$. We apply these established results in precisely the following setting.

\begin{proposition}\label[proposition]{prop:spectral}
For every prime $p\equiv1\pmod {17}$ and every $n\geq1$, each eigenvalue $\lambda$ of $A_{n,p}$ with an eigenvector in $\one^\perp$ satisfies $|\lambda|\leq2\sqrt p$.
\end{proposition}
\begin{proof}
For functions on ideal classes, the arithmetic inner product uses the reciprocal of the order of each projective unit stabilizer as its weight. All these stabilizers are trivial by \cref{prop:class}, so it is the counting inner product defined above. Moreover, the reverse-edge pairing in \cref{eq:rotation} makes $A_{n,p}$ symmetric. The identification in \cref{prop:hecke} therefore makes $T_P$ self-adjoint for this inner product.

The quaternionic form of Jacquet--Langlands separates off the one-dimensional constituents given by \emph{norm characters}, namely functions $h([g])=\chi((\nrd(g_R))_R)$ for a multiplicative character~$\chi$, see \cite[Def.~3.7 and Thm~3.9]{DV13} and \cite[\S7]{KiV10}. This classification is part of the external automorphic input. We now show directly that every such function in our class space is constant.

At every finite prime ideal $R$, the reduced norms of $(\EE_n)_R^\times$ are all of $\mathcal O_{F,R}^\times$. Indeed, after conjugation the Eichler order contains $\diag(u,1)$ and its inverse for every local scalar unit $u$, and reduced norm is determinant. Hence
\[
 \nrd(\mathcal U_n)=\prod_R\mathcal O_{F,R}^\times
                  =\nrd(\mathcal U_0).
\]
A norm-character function invariant under $\mathcal U_n$ is therefore invariant under $\mathcal U_0$. It descends to
\[
 D^\times\backslash\widehat D^\times/\mathcal U_0
 \simeq\Cl(\EE_0),
\]
which is a singleton by \cref{prop:class}. Thus the noncuspidal subspace consists exactly of constants. This argument does not require a description of the global reduced-norm image.

The central action is trivial as well. Since the ordinary ideal class group of $F$ is trivial, each finite scalar idele $a=(a_R)_R$ can be written $a=bu$ with $b\in F^\times$ and $u\in\prod_R\mathcal O_{F,R}^\times$. Multiplication by $b$ is absorbed on the left by $D^\times$, and the scalar unit tuple $u$ belongs to every $\mathcal U_n$. Thus scalar ideles act trivially on the class set.

By Jacquet--Langlands, $\one^\perp$ with its Hecke action consequently embeds in the space of Hilbert cusp forms of parallel weight $(2,2)$, level $Q^{n+2}$, and trivial central character; see \cite[Thm~3.9]{DV13}. Here ``embeds with its Hecke action'' means that the linear embedding commutes with each Hecke operator away from the level. Since $P\ne Q$, there is no level condition at $P$. In representation-theoretic terminology, the local representation (a linear action of the local matrix group) is \emph{unramified}: it has a nonzero vector invariant under the full integral matrix group $\GL_2(\mathcal O_{F,P})$.

Blasius's Hilbert modular Ramanujan theorem \cite[Thm~1]{Bla06} gives
\[
 \beta+\gamma=\lambda_P,\qquad
 \beta\gamma=p,\qquad |\beta|=|\gamma|=\sqrt p.
\]
Hence $|\lambda_P|\leq2\sqrt p$. The Hecke operators are self-adjoint, meaning equal to their adjoints for the stated inner product, and the bound holds on the full complement of constants. The identification of adjacency with $T_P$ in \cref{prop:hecke} proves the claim.
\end{proof}

Since $p+1>2\sqrt p$, the eigenvalue $d=p+1$ has multiplicity one: any additional eigenvector with eigenvalue $d$ could be chosen orthogonal to the constants, contradicting the bound. This proves connectedness, and the absence of $-d$ proves non-bipartiteness. Dividing adjacency by $d$ gives the random-walk operator; its largest absolute nonconstant eigenvalue is at most $2\sqrt{d-1}/d$. Its \emph{absolute spectral gap}, one minus that quantity, is therefore at least
\begin{equation}\label{eq:gap}
 1-\frac{2\sqrt{d-1}}d.
\end{equation}
This proves parts (i) and (ii) of \cref{thm:main}.
This completes Step~3 of the proof outlined in Section~\ref{sec:outline}

\bigskip\section{Binary covers and an explicit signing}\label{sec:lifts}

The spectral estimate does not by itself show that consecutive graphs form a $2$-lift. We prove that each neighbor transformation is an \emph{isometry} for the $2$-adic distance $|t-u|_2=2^{-v_2(t-u)}$, defined in $\S$\ref{subsec:completion}: it preserves that distance exactly. This makes reduction a label-preserving two-sheeted cover, and the next output bit gives the straight/crossed signing without any search.

\begin{lemma}[{\rm Dyadic isometry}{}]\label[lemma]{lem:isometry}
For every $g=\left(\begin{smallmatrix}a&b\\c&e\end{smallmatrix}\right)$ arising in \cref{eq:gj}, and every distinct $t,u\in\Z_2$,
\begin{equation}\label{eq:isometry}
 v_2(\mobius_g(t)-\mobius_g(u))=v_2(t-u).
\end{equation}
\end{lemma}
\begin{proof}
Direct subtraction gives
\[
 \mobius_g(t)-\mobius_g(u)
 =\frac{(ae-bc)(t-u)}{(2ct+e)(2cu+e)}.
\]
The determinant and both denominator factors are odd. This proves the valuation identity. The inverse matrix has the same scalar-modulo-$2$ property, so $\mobius_g$ is a bijection of $\Z_2$.
\end{proof}

The isometry property gives more than compatibility with reduction: the transformation for each edge label maps the two children of a vertex bijectively onto the two children of its neighbor. The pairing of reverse-edge labels is independent of the level. Therefore reduction is a two-sheeted label-preserving cover in the graph convention of \cref{s:notation}.

The signing can be written without any search. Choose the representative $0\leq t<2^n$, and define
\begin{align}
 y&=f_j(t)\bmod2^n,\label{eq:liftbit1}\\
 \beta_n(t,j)&=\frac{(f_j(t)\bmod2^{n+1})-y}{2^n}\in\{0,1\},
 \qquad \sigma_n(t,j)=(-1)^{\beta_n(t,j)}.\label{eq:liftbit}
\end{align}
The notation $\bmod$ here denotes the least nonnegative residue. The two-sheet labeling is $t+b2^n\leftrightarrow(t,b)$, with $b\in\{0,1\}$. The bit $\beta_n(t,j)=0$ means a straight lift, while $1$ means a crossed lift; $\sigma_n(t,j)\in\{+1,-1\}$ is the corresponding sign. We use $b\mathbin\oplus\beta$ for addition modulo $2$.

\begin{proposition}[{\rm Straight/crossed rule}{}]\label[proposition]{prop:lift}
With the two-sheet labeling, the rotation map at the next level is
\begin{equation}\label{eq:liftformula}
 \operatorname{Rot}_{n+1,p}\bigl((t,b),j\bigr)
 =\bigl((y,b\mathbin\oplus\beta_n(t,j)),j^*\bigr).
\end{equation}
Moreover, $\beta_n(t,j)=\beta_n(y,j^*)$, so the signs agree on the two orientations of each edge. Define the \emph{signed adjacency operator} by
\[
 (B_{n,p}h)(t)=\sum_{j=1}^{p+1}\sigma_n(t,j)h(f_j(t)\bmod2^n).
\]
Then
\begin{equation}\label{eq:spectraldecomp}
 A_{n+1,p}\simeq A_{n,p}\oplus B_{n,p},
 \qquad \|B_{n,p}\|\leq2\sqrt p.
\end{equation}
Here $\oplus$ is a direct sum of linear operators, $\simeq$ means equivalence after an orthogonal change of basis, and $\|B_{n,p}\|$ is its operator norm for the counting inner product, equivalently its largest absolute eigenvalue.
\end{proposition}
\begin{proof}
The difference identity in \cref{lem:isometry}, with $u=t+2^n$, gives
$f_j(t+2^n)-f_j(t)\equiv2^n\pmod {2^{n+1}}$ because the quotient multiplying $2^n$ is odd. This proves \cref{eq:liftformula}. Applying the inverse permutation to the lifted edge proves the reverse-edge identity.

The functions on the lifted vertices split into functions constant on each fiber and functions that change sign between its two points. In the first subspace the adjacency action is $A_{n,p}$. In the second it is $B_{n,p}$. This proves the direct sum. Every fiber-antisymmetric function is orthogonal to the global constant function, so \cref{prop:spectral} at level $n+1$ gives the bound on $B_{n,p}$.
\end{proof}

For edges between distinct base vertices this is literally the Bilu--Linial parallel/crossed matching rule. A half-loop can lift to two half-loops or to an edge joining the two sheets. This explains both the precise relationship with $2$-lifts and the graph-convention qualification. Part (iii) of \cref{thm:main} follows.
This completes Step~4 of the proof outlined in Section~\ref{sec:outline}

\bigskip\section{Strong explicitness and bit complexity}\label{sec:complexity}

To establish strong explicitness, we must compute the fixed transformations and evaluate them at a prescribed vertex without constructing the graph. We separate the degree-dependent preprocessing from the level-dependent query, and bound both in the bit model. A fixed norm-lattice certificate makes the finite search bounds explicit.

\subsection{Deterministic preprocessing in the degree}
The preprocessing searches the fixed rank-eight integer lattice underlying \cref{eq:E0basis}. The search depends on $p$, but does not enumerate vertices of any graph.

\begin{lemma}[{\rm A balanced generator}{}]\label[lemma]{lem:balanced}
A totally positive generator $\pi$ of $\adjprime$ can be chosen so that $\pi_1=\tau_+(\pi)$ and $\pi_2=\tau_-(\pi)$ obey
\begin{equation}\label{eq:balanced}
 \fundunit^{-1}\sqrt p\leq\pi_i\leq\fundunit\sqrt p,
 \qquad \fundunit=4+\sqrt{17}.
\end{equation}
It can be found deterministically in $O(p)\thinspace \textrm{\em poly}(\log p)$ bit operations from $p$ and the residue $s_p$ in \cref{eq:activeprime}.
\end{lemma}
\begin{proof}
Start with any totally positive generator and multiply by a power of the totally positive unit~$\fundunit^2$. The first embedding is multiplied by $\fundunit^{2k}$ and the second by $\fundunit^{-2k}$. A choice of $k$ places the first embedding in the interval in \cref{eq:balanced}. The product of the embeddings is $p$, so the second is then in the same interval.

This also gives a finite search not requiring prior principalization. Write $\pi=a+b\omega$. Search $|a|,|b|\leq20\lceil\sqrt p\rceil$ and test
\begin{equation}\label{eq:generatorsearch}
 a^2+ab-4b^2=p,\quad 2a+b>0,\quad
 (2a+b)^2\leq68p,\quad a+b s_p\equiv0\pmod p.
\end{equation}
Balanced generators occur in this box, by solving for $a,b$ from their two embeddings. The first two tests imply total positivity. To interpret the trace bound, put $x=\pi_1/\sqrt p>0$, so $\pi_2/\sqrt p=x^{-1}$. Since $\fundunit+\fundunit^{-1}=2\sqrt{17}$, we have
\[
 (\pi_1+\pi_2)^2\leq68p
 \quad\Longleftrightarrow\quad x+x^{-1}\leq2\sqrt{17}
 \quad\Longleftrightarrow\quad\fundunit^{-1}\leq x\leq\fundunit.
\]
Thus the trace test is exactly the required bound on the two individual embeddings; their ratio lies in $[\fundunit^{-2},\fundunit^2]$.

Finally, containment in $\adjprime$ and equality of norms imply $(\pi)=\adjprime$. The box has $O(p)$ points, and every test is exact arithmetic on $O(\log p)$-bit integers. Choose the lexicographically first solution to fix a deterministic convention.
\end{proof}

Let the fixed $\Z$-basis be
\begin{equation}\label{eq:integerbasis}
 \mathcal B=(1,\omega,\eta\qi,\omega\eta\qi,
                  \eta u,\omega\eta u,v,\omega v).
\end{equation}
For $x=\sum_{s=1}^8 m_s\mathcal B_s$, let $M$ be the Gram matrix of the real quadratic form obtained by taking the field trace of the reduced norm. Thus $M$ is the rational positive-definite matrix satisfying
\[
 \Tr_{F/\Q}\nrd(x)=m^{\mathsf T}M m.
\]
An exact matrix certificate is given at the end of this section, in \cref{app:gram}. Its inverse has maximal diagonal entry $828/17$. If $\nrd(x)=\pi$ with $\pi$ balanced, Cauchy--Schwarz for the $M$-inner product yields
\begin{align}
 m_s^2
 &\leq (M^{-1})_{ss}\,m^{\mathsf T}M m\notag \leq \frac{828}{17}\,2\sqrt{17}\sqrt p
 <21^2\sqrt p.\label{eq:boxbound}
\end{align}
Thus every solution is found in the completely specified box
\begin{equation}\label{eq:enumerationbox}
 |m_s|\leq\lceil21p^{1/4}\rceil\qquad(1\leq s\leq8).
\end{equation}

\begin{proposition}[{\rm Computing the generators}{}]\label[proposition]{prop:compiler}
Enumerating \cref{eq:enumerationbox}, testing the exact equation $\nrd(x)=\pi$, and retaining the sign for which the first nonzero $m_s$ is positive produces $d=p+1$ representatives for $S_p$. The running time is $O(p^2\operatorname{poly}(\log p))$. Each generator has $O(\log p)$ bits of coefficient data in the fixed basis. Sorting the surviving coordinate vectors lexicographically fixes the edge-label ordering. Conjugating each vector, normalizing its sign by the same rule, and locating it in the sorted list computes the reverse-edge label $j^*$.
\end{proposition}
\begin{proof}
The box contains $O(p^{8/4})=O(p^2)$ integer points, with an absolute constant because the field and the order are fixed. Reduced norm is an explicit quadratic form over $F$. Each test takes $\operatorname{poly}(\log p)$ bit operations. Bound \cref{eq:boxbound} proves completeness, and \cref{lem:normgenerators} proves that exactly $d$ sign classes survive. The coordinate bound also gives the coefficient-height assertion. Quaternion conjugation preserves $\EE_0$ and is represented by a fixed integral matrix in the basis~$\mathcal B$. Applying that matrix, normalizing the sign, and using binary search in the sorted list takes $O(p\operatorname{poly}(\log p))$ additional bit operations over all labels. This fits within the stated preprocessing bound and gives the complete reverse-label table.
\end{proof}

The search box has large constants, but gives an explicit stopping rule and a uniform bound in $p$. It avoids a separate principalization algorithm.

\subsection{The local adjacency oracle}
In the fixed basis $\mathcal B=(\mathcal B_1,\ldots,\mathcal B_8)$ of \cref{eq:integerbasis}, write each quaternion as
\[
 x_j=\sum_{s=1}^{8}m_{js}\mathcal B_s,
 \qquad (m_{j1},\ldots,m_{j8})\in\Z^8.
\]
These eight integers specify $x_j$ exactly. Together with the fixed matrix representation $\splitmap$, they determine $g_j=\splitmap(\bar x_j)$ and the neighbor function $f_j$ at every level. Store the coordinate vectors and reverse-label table once, independently of $n$; no infinite $2$-adic expansion is stored.

A query is specified by $n,t,j$ and whether it requests a neighbor or a lift bit. Set the target output precision to $k=n$ for a neighbor query and $k=n+1$ for a lift-bit query. Compute the selected matrix $g_j$ modulo $2^{k+1}$, and evaluate
\begin{equation}\label{eq:oracle}
 u=(a_jt+b_j/2)(2c_jt+e_j)^{-1}\pmod {2^k},
 \qquad 0\leq u<2^k.
\end{equation}
Here $b_j$ is represented by its even residue modulo $2^{k+1}$, so $b_j/2$ is well-defined modulo $2^k$. The denominator is odd and hence invertible. For a neighbor query return $u$. For a lift query return $\beta_n(t,j)=\lfloor u/2^n\rfloor$ and, when requested, the sign $(-1)^{\beta_n(t,j)}$. In particular, a lift query needs matrix entries modulo $2^{n+2}$, not merely modulo $2^{n+1}$.

The arithmetic procedure is as follows.
\begin{enumerate}[label=\textup{\arabic*.}]
\item Set $k$ as above, and lift the even roots $\alpha,z$ of $X^2+3X-2$ and $X^2+X+2$ to precision $2^{k+1}$. Put $\omega=\alpha+2$ and $\eta=\alpha+4$ at this precision.
\item Reconstruct the matrix of $x_j$ from its eight integer coordinates using the basis in \cref{eq:integerbasis} and the integral formulas for $\splitmap(\qi),\splitmap(u),\splitmap(v)$ in \cref{prop:explicit-order}. These formulas have only the odd denominator $3$, which is inverted modulo $2^{k+1}$. Take the matrix adjugate to obtain $g_j=\splitmap(\bar x_j)$ modulo $2^{k+1}$.
\item Divide the even upper-right entry by $2$, and compute the numerator and odd denominator in \cref{eq:oracle} modulo $2^k$. Use modular inversion to obtain $u$.
\item Return the neighbor residue or extract its $n$th binary digit as specified above. The reverse-edge label is read from the precomputed table.
\end{enumerate}
The additional matrix bit is necessary because of the division by $2$ in step 3. Using the integral basis-image formulas in step 2 avoids any further loss of binary precision. An implementation that instead substitutes the original rational quaternion expressions must retain the fixed additional guard bits required before divisions by their bounded powers of $2$.

\begin{proposition}[{\rm Local complexity}{}]\label[proposition]{prop:local}
Given the integer coordinates of the selected quaternion, the preceding algorithm computes a neighbor or a lift bit in $O((n+\log p)^3)$ bit operations, with polynomial bit space. For fixed $p$ this is $O_p(n^3)$, with one finite algorithm valid for all levels.
\end{proposition}
\begin{proof}
Put $L=n+\log p+O(1)$. A simple root modulo $2^s$ extends uniquely to one of $x$ and $x+2^s$ modulo $2^{s+1}$. Testing these two possibilities for each successive bit uses $O(L)$ arithmetic steps. With schoolbook multiplication, each step costs $O(L^2)$ bit operations. Both roots are therefore available in $O(L^3)$ time.

Evaluating the fixed-size matrix expressions and the M\"obius formula uses a bounded number of operations on $O(L)$-bit integers. The extended Euclidean algorithm, or Newton inversion, computes the odd modular inverse within the same cubic upper bound. Formula \cref{eq:liftbit} changes the requested precision by only one bit. At no point is an adjacency matrix, a list of earlier levels, or a signing table of size $2^n$ constructed.
\end{proof}

The degree description has $O(d\log p)$ bits. If stored as an indexed read-only table, a query reads only the eight integer coordinates of one quaternion. On a sequential model one may instead absorb the degree-dependent lookup cost into the constant in $O_p(n^3)$. The statement of strong explicitness for fixed degree is unchanged.

The preprocessing can even be rerun on every query and still be a constant cost as a function of $n$ for a fixed admissible degree. Thus the construction does not rely on a nonuniform sequence of good seeds. Conversely, our bounds do \emph{not} give a single neighbor algorithm with no degree-dependent preprocessing and running time polynomial in $(\log d,\log N)$. This distinction is part of the theorem, not an omitted implementation detail. Part (iv) of \cref{thm:main} now follows.

\subsection{An exact certificate for the norm-enumeration bound}\label{app:gram}

This subsection supplies a finite certificate for the enumeration bound. Use the basis \cref{eq:integerbasis} and define $H=2M$. Direct expansion of \cref{eq:norm} gives
\[
H=\begin{pmatrix}
4&2&0&0&5&11&0&-17\\
2&18&0&0&11&31&-17&-17\\
0&0&42&106&0&0&5&11\\
0&0&106&274&0&0&11&31\\
5&11&0&0&42&106&0&0\\
11&31&0&0&106&274&0&0\\
0&-17&5&11&0&0&24&12\\
-17&-17&11&31&0&0&12&108
\end{pmatrix}.
\]
In particular $\Tr\nrd(x)=\tfrac12m^{\mathsf T}Hm$. The diagonal entries of $M^{-1}=2H^{-1}$ are
\[
 \frac1{17}(828,184,137,21,822,126,108,24).
\]
The form is positive definite because it is the sum of the two positive-definite real norm forms. Thus Cauchy--Schwarz applies without any numerical eigenvalue estimate. Since
\[
 \frac{828}{17}\,2\sqrt{17}<441,
\]
the coefficient bound $21p^{1/4}$ used in \cref{eq:enumerationbox} is rigorous.
This completes Step~5 of the proof outlined in Section~\ref{sec:outline}.

\bigskip\section{Consequences of the binary organization}\label{sec:consequences}

The binary organization makes random walks and eigenfunctions compatible across all levels. We record these consequences, the complete-graph sparsifier estimate, and the limits of the construction. The number-theoretic comparison with LPS is given in~$\S$\ref{subsec:lps-comparison}. Throughout this section, prime $p$ is fixed and satisfies $p\equiv1\pmod {17}$.

\subsection{Exactly compatible random walks}
Let $W_n=A_{n,p}/(p+1)$ be the random-walk operator. Write $\red_n:V_{n+1}\to V_n$ for reduction modulo $2^n$, and define the pullback $J_nh=h\circ\red_n$. Each labeled neighbor transformation respects reduction, so
\begin{equation}\label{eq:walkprojection}
 W_{n+1}J_n=J_nW_n.
\end{equation}
Indeed, both sides at a vertex $t$ equal the average of $h$ over the projected neighbors of $t$, with the same edge-incidence multiplicities. Thus a walk generated at level $n+1$ projects to exactly the level-$n$ walk, not merely to a process with approximately the same distribution. Sharing the edge-label choices gives an explicit coupling across all levels. A uniformly random fine-level vertex has a uniformly random reduction, since every fiber has size two.

This is the additional structure beyond an unrelated sequence of LPS or Morgenstern graphs. It can be useful when an algorithm maintains a state at several binary resolutions. A separate application analysis would still be required to turn that compatibility into better randomness, time, or space bounds.

\subsection{Old eigenvalues never disappear}
Write $\lambda_*(G)$ for the largest absolute adjacency eigenvalue on the mean-zero functions, and call a function mean-zero when the sum of its values is zero. If $A_{n,p}h=\lambda h$, the covering identity gives $A_{n+1,p}J_nh=\lambda J_nh$. A mean-zero function remains mean-zero under pullback. Consequently, a nontrivial eigenvalue introduced at one level persists at every higher level. In any fixed-degree tower with this property, full-spectrum estimates of the form
\[
 \lambda_*(G_n)\leq2\sqrt p+\varepsilon_n,\qquad \varepsilon_n\longrightarrow0,
\]
would already force every level to obey the exact bound: an excess $\delta>0$ at one level would imply $\varepsilon_n\geq\delta$ at all later levels. This argument does not apply to unrelated graphs, a fixed positive error, or a bound on only the newly introduced eigenvalues.

\subsection{Locally accessible complete-graph sparsifiers}
Put $N=2^n$, $d=p+1$, and let $L_n=dI-A_{n,p}$ be the graph Laplacian, with $I$ the identity matrix. Diagonal adjacency contributions cancel in $L_n$, as they should for loops and half-loops. On $\one^\perp$ the eigenvalues of $L_n$ are in $[d-2\sqrt p,d+2\sqrt p]$. Since the Laplacian of the complete graph $K_N$ is $L_{K_N}=NI-J$, with $J$ the all-ones matrix, it follows that
\begin{equation}\label{eq:sparsifier}
 (1-\varepsilon_p)L_{K_N}
 \preceq \frac{N}{p+1}L_n
 \preceq (1+\varepsilon_p)L_{K_N},
 \qquad \varepsilon_p=\frac{2\sqrt p}{p+1}.
\end{equation}
The notation $X\preceq Y$ means that $v^{\mathsf T}(Y-X)v\geq0$ for every real vector $v$ (positive-semidefinite order). The inequalities hold on the full space, because both Laplacians annihilate constants. After discarding diagonal loops and combining parallel edges, the matrix describes a weighted graph with $O(pN)$ edges and locally computable edge incidences. This is the standard Ramanujan sparsifier example \cite{BSS12}, now available at each binary size. 
To enumerate a whole row, query all $p+1$ edge labels and aggregate their contributions.


\vskip.5cm
{\small 
\noindent\textbf{AI methodology.}
\noindent
OpenAI's ChatGPT 6 was substantially used as an editing tool. It helped to bridge the distance between the mathematical and computer science presentations.
Notably, the original approach in [KV26a, MP17] was
modified to make the paper self-contained albeit heavily number-theoretic, at the expense of
losing the motivating flavor of automata groups.  Additionally, ChatGPT was used as a library assistant to assemble, format and organise the references.
The authors take full responsibility for the correctness of all results.
}
\vskip.7cm


{

}

\vskip.4cm

\end{document}